\documentclass{article}
\usepackage {amsmath, amsfonts, amssymb, mathrsfs, amsthm, diagrams,accents, sectsty,hyphenat}
\usepackage[usenames]{color}
\usepackage{palatino}
\def\mf#1{\mathfrak{#1}}
\def\mc#1{\mathcal{#1}}
\def\mb#1{\mathbb{#1}}
\def\tx#1{{\rm #1}}
\def\tb#1{\textbf{#1}}
\def\ti#1{\textit{#1}}
\def\tr{\tx{tr}\,}
\def\R{\mathbb{R}}
\def\C{\mathbb{C}}
\def\Q{\mathbb{Q}}

\def\Z{\mathbb{Z}}

\def\ol#1{\overline{#1}}
\def\ul#1{\underline{#1}}

\def\hat{\widehat}
\def\vol{\tx{vol}}

\newarrow{Equals}{=}{=}{=}{=}{}

\def\rmk{\tb{Remark: }}

\def\pf{\tb{Proof: }}

\newenvironment{mytitle}
{\begin{center}\large\sc}
{\end{center}}

\newtheorem{thm}{Theorem}[subsection]
\newtheorem{lem}[thm]{Lemma}
\newtheorem{pro}[thm]{Proposition}
\newtheorem{cor}[thm]{Corollary}

\newtheorem{fct}[thm]{Fact}

\sectionfont{\center\sc\normalsize}
\subsectionfont{\bf\normalsize}

\numberwithin{equation}{section}

\def\phi{\varphi}

\begin{document}

\begin{mytitle} Endoscopic character identities for depth-zero supercuspidal $L$-packets \end{mytitle}
\begin{center} Tasho Kaletha \end{center}

In this paper we prove the conjectural endoscopic character identities for the local Langlands correspondence constructed in \cite{DR09}. The local Langlands
correspondence, which is known in the real case and partially constructed in the p-adic case, assigns to each Langlands parameter for a reductive group $G$
over a local field $F$ a finite set of admissible irreducible representations of $G(F)$, called an $L$-packet. When such a parameter factors through an
endoscopic group $H$, the broad principle of Langlands functoriality asserts that the packet on $H$ should "transfer" to the packet on $G$. The endoscopic
character identities are an instance of this principle -- they state that the "stable" character of the packet on $H$ is identified via endoscopic induction
with an "unstable" character of the packet on $G$.

To be more precise, let $F$ be a p-adic field with Weil-group $W_F$ and let $G$ be a connected reductive group over $F$. For the purposes of this introduction,
we assume that $G$ is unramified, although in the body of this paper the more general case of a pure inner form of an unramified group is handled. Let ${^LG}$
be an $L$-group for $G$, that is ${^LG} = \hat G \rtimes W_F$, where $\hat G$ is the complex Langlands dual of $G$ and $W_F$ acts on $\hat G$ via its action on
the based root datum of $\hat G$ which is dual to that of $G$. The Langlands parameters considered in this paper are continuous sections
\[W_F \rightarrow {^LG} \]
of the natural projection ${^LG} \rightarrow W_F$ and subject to certain conditions, called TRSELP in \cite{DR09}, which will be reviewed in detail later on.
To such a parameter DeBacker and Reeder construct in loc.cit. an $L$-packet $\Pi_G(\phi)$ of representations of $G(F)$ and a bijection
\[ \tx{Irr}(C_\phi,1) \rightarrow \Pi_G(\phi), \qquad \rho \mapsto \pi_\rho \]
where $C_\phi$ is the component group of the centralizer in $\hat G$ of $\phi$ and $\tx{Irr}(C_\phi,1)$ are those representations of the finite group $C_\phi$
which are trivial on elements of $C_\phi$ coming from the center of $\hat G$. This bijection maps the trivial representation of $C_\phi$ to a generic
representation of $G(F)$.

Let $(H,s,\hat\eta)$ be an unramified endoscopic triple for $G$. Recall that $H$ is an unramified reductive group over $F$, $s$ is a Galois-fixed element of
the center of $\hat H$, and $\hat\eta$ is an inclusion $\hat H \rightarrow \hat G$ which identifies $\hat H$ with $(\hat G_{\hat\eta(s)})^\circ$. It was shown
by Hales that $\hat\eta$ extends to an embedding ${^L\eta} : {^LH} \rightarrow {^LG}$. Thus for any parameter $\phi^H$ for $H$ we may consider the parameter
$\phi = {^L\eta} \circ \phi^H$, i.e. we have
\begin{diagram}
^LH&\rTo^{^L\eta}&^LG\\
\uTo<{\phi^H}&\ruTo>{\phi}\\
W_F
\end{diagram}
If both parameters are of the type considered here, then we have the $L$-packets $\Pi_G(\phi)$ and $\Pi_H(\phi^H)$. Associated to these, we have the stable
character
\[ \mc{S}\Theta_{\phi^H} := \sum_{\rho \in \tx{Irr}(C_{\phi^H},1)} [\dim\rho] \chi_{\pi_\rho} \]
of $\Pi_H(\phi^H)$, which is a stable function on $H(F)$ (this is one of the main results of \cite{DR09}), as well as the $s$\-unstable character
\[ \Theta^s_{\phi,1} := \sum_{\rho \in \tx{Irr}(C_\phi,1)} [\tr \rho(s)] \chi_{\pi_\rho} \]
of $\Pi_G(\phi)$, which is an invariant function on $G(F)$.

Recall that the representation $\pi_1$ of $G(F)$ is generic. Thus there is a Borel subgroup $B=TU$ of $G$ defined over $F$ and a generic character $\psi : U(F)
\rightarrow \C^\times$ which occurs in the restriction of $\pi_1$ to $U(F)$. Associated to the character $\psi$ there is a unique normalization $\Delta_\psi$
of the transfer factor for $G$ and $H$, called the Whittaker normalization. The endoscopic lift of the stable function $\mc{S}\Theta_{\phi^H}$ is given by
\[ \tx{Lift}^G_H\mc{S}\Theta_{\phi^H}(\gamma) := \sum_{\gamma^H} \Delta_\psi(\gamma^H,\gamma)\frac{D(\gamma^H)^2}{D(\gamma)^2}\mc{S}\Theta_{\phi^H}(\gamma^H) \]
where $\gamma \in G(F)$ is any strongly regular semi-simple element and $\gamma^H$ runs through the set of stable classes of $G$-strongly regular semi-simple
elements in $H(F)$.

The main result of this paper asserts that
\[ \Theta^s_{\phi,1} = \tx{Lift}^G_H\mc{S}\Theta_{\phi^H} \]
As a corollary of the main result in the case where $G$ is a pure inner form of an unramified group $G^*$ and $H=G^*$ we obtain a proof (for the $L$-packets
considered) of the conjecture of Kottwitz \cite{Kot83} about sign changes in stable characters on inner forms.

We now describe the contents of the paper. After fixing some basic notation in Section \ref{sec:not}, we discuss pure inner twists and the associated notions
of conjugacy and stable conjugacy. We have allowed trivial inner twists in the discussion so as to accommodate the natural construction of the $L$-packets in
\cite{DR09} and not just their normalized form. With these notions in place we implement an observation of Kottwitz which allows one to define compatible
normalizations of the absolute transfer factors for all pure inner twists. In Section \ref{sec:main} we briefly review the construction of the local Langlands
correspondence in \cite{DR09}, and after gathering the necessary notation we state the main result of this paper. The remaining sections are devoted to its
proof, which is similar in spirit to the proof of the stability result in loc. cit.. In Section \ref{sec:signs} we study three signs which are defined for a
pair $(G,H)$ of a group $G$ and an endoscopic group $H$ and play an important role in the theory of endoscopy -- one of them is defined in terms of the split
ranks of these groups and goes back to \cite{Kot83}, the other one occurs in Waldspurger's work \cite{Wal95} on the endoscopic transfer for p-adic Lie
algebras, and the third is a certain local $\epsilon$-factor used in the Whittaker normalization of the transfer factors \cite{KS99}. We show that when both
$G$ and $H$ are unramified, these three signs coincide. This supplements the results of \cite[\S12]{DR09} to assert in particular that the Waldspurger-sign and
the relative-ranks-sign coincide whenever $G$ is a pure inner form of an unramified group and $H$ is an unramified endoscopic group. Because this section may
be of independent interest we have minimized the notation that it borrows from previous sections. Section \ref{sec:schar} deals with establishing a reduction
formula for the unstable character of an $L$-packet with respect to the topological Jordan decomposition. For that we first need explicit formulas for some
basic constructions in endoscopy, which are established in two preparatory subsections. Among other things we show that the isomorphism $H^1(F,G) \rightarrow
\tx{Irr}(\pi_0(Z(\hat G)^\Gamma))$ constructed in \cite{DR09} via Bruhat-Tits theory coincides with the one constructed in \cite{Kot86} using Tate-Nakayama
duality. With these preliminaries in place we derive the reduction formula for the unstable character using the results of \cite[\S9,\S10]{DR09}. The
ingredients from the previous sections are combined in Section \ref{sec:chid} to establish the proof of the main result. After reducing to the case of compact
elements the reduction formula from Section \ref{sec:schar} is combined with the work of Langlands and Shelstad \cite{LS90} and Hales \cite{Hal93} on
endoscopic descent. The topologically unipotent part of the resulting expression is then transferred to the Lie algebra, where we invoke the deep results of
Waldspurger on endoscopic transfer for p-adic Lie-algebras together with the fundamental lemma, which has been recently proved by the combined effort of many
people.

We would like to bring to the attention of the reader some related work on this problem. In \cite{KV1}, Kazhdan and Varshavsky construct an endoscopic
decomposition for the $L$-packets considered here. In particular, they consider the $s$-unstable characters of these packets and show that they belong to a
space of functions which contains the image of endoscopic induction. The existence of such a decomposition is a necessary condition for the validity of the
character identities considered here and also gave us yet more reason to hope that indeed these identities should be true. In \cite{KV2} the aforementioned
authors prove a formula for the geometric endoscopic transfer of Deligne-Lusztig functions, in particular answering a conjecture of Kottwitz. After the current
paper was written, the author was informed in a private conversation with Kazhdan that the results in \cite{KV2} could likely be used to derive character
identities similar to the ones proved here, at least on the set of elliptic elements, and possibly in general.

The author would like to thank Professor Robert Kottwitz for his generous support and countless enlightening and inspiring discussions. This work would not
have been possible without his dedication and kindness. The author would also like to thank Professor Stephen DeBacker for suggesting this problem and
discussing at length the constructions and character formulas in \cite{DR09}, as well as for his continual support and encouragement.

\newpage
\tableofcontents
\newpage


\section{Notation} \label{sec:not}

Let $F$ be a $p$-adic field (i.e. a finite extension of $\Q_p$) with ring of integers $O_F$, uniformizer $\pi_F$, and residue field $k_F=O_F/\pi_F O_F$ with
cardinality $q_F$. We use analogous notation for any other discretely valued field, in particular for the maximal unramified extension $F^u$ of $F$ in a fixed
algebraic closure $\ol{F}$. Since we will consider only extensions of $F$ which lie in $F^u$, $\pi_F$ will be a uniformizer in each of them and so we will drop
the index $F$ and simply call it $\pi$. For any such finite extension $E$, $v_E : E^\times \rightarrow \Z$ will be the discrete valuation normalized so that
$v_E(\pi)=1$, and $|x|_E$ will be the norm given by $q_E^{-v_E(x)}$. Thus $v_E$ extends $v_F$ and so we may again drop the index $F$. On the other hand, for $x
\in F^\times$ we have $|x|_E=|x|_F^{[E:F]}$; if $dx$ is any additive Haar measure on $E$ then $d(ax)=|a|_Edx$. The absolute Galois group of $F$ will be denoted
by $\Gamma$, its Weil group by $W_F$ and inertia group by $I_F$. We choose an element $\tx{Fi} \in \Gamma$ whose inverse induces on $\ol{k_F}$ the map $x
\mapsto x^{q_F}$.

For a reductive group $G$ defined over $F$, we will denote its Lie algebra by the Fraktur letter $\mf{g}$. Our convention will be that $a \in G$ resp. $a \in
\mf{g}$ will mean that $a$ is an $\ol{F}$-point of the corresponding space, while a maximal torus $T \subset G$ will be tacitly assumed to be defined over $F$.
The action of $\tx{Fi}$ on both $G(F^u)$ and $\mf{g}(F^u)$ will be denoted by $\tx{Fi}_G$. For a semi-simple $a \in G$, we will write $\tx{Cent}(a,G)=G^a$ for
the centralizer of $a$ in $G$ and $G_a$ for its connected component. If $T \subset G$ is a maximal torus then the roots resp. coroots of $T$ in $G$ will be
denoted by $R(T,G)$ resp. $R^\vee(T,G)$. The center of $G$ will be $Z_G$, or simply $Z$ if $G$ is understood, and the maximal split torus in $Z_G$ will be
$A_G$. The sets of strongly-regular semi-simple elements of $G$ resp. $\mf{g}$ will be denoted by $G_\tx{sr}$ resp. $\mf{g}_\tx{sr}$. The set of compact
elements in $G(F)$ will be denoted by $G(F)_0$ (note that we are using the wording of \cite{DR09} here; in \cite{Hal93} these elements are called
strongly-compact). For any $g \in G$ the map $G \rightarrow G, x \mapsto gxg^{-1}$ as well as its tangent map $\mf{g} \rightarrow \mf{g}$ will be called
$\tx{Ad}(g)$. Abusing words, will will refer to the orbits of $\tx{Ad}(G)$ in $\mf{g}$ as conjugacy classes, and then notions such as stable classes and
rational classes will have their obvious meaning.

To maintain notational similarity with \cite{DR09}, we will sometimes use the following conventions. If $\psi : G \rightarrow G'$ is an inner twist, then we
may identify $G(\ol{F})$ and $G'(\ol{F})$ via $\psi$ and suppress $\psi$ from the notation, thereby treating $\gamma \in G(\ol{F})$ and $\psi(\gamma) \in
G'(\ol{F})$ as the same element. If $u \in Z^1(\Gamma,G)$ is a cocycle, then we will use the same letter $u$ also for the value of that cocycle at $\tx{Fi}$.

If $(H,s,\hat\eta)$ is an endoscopic triple for a reductive group $G/F$, we will often attach a superscript $H$ to objects related to $H$, such as maximal
tori, Borels, or elements of $H(F)$. If $^L\eta : {^LH} \rightarrow {^LG}$ is an $L$-embedding extending $\hat\eta$, then we will call $(H,s,{^L\eta})$ an
extended triple for $G$. The set of $G$-strongly regular semi-simple elements of $H$ resp. $\mf{h}$ will be denoted by $H_{G-\tx{sr}}$ resp.
$\mf{h}_{G-\tx{sr}}$. Let $t^H \in H(F)$ and $t \in G(F)$ be semi-simple elements. We will call $t$ an image of $t^H$ if there exist maximal tori $T^H \subset
H$ and $T \subset G$ and an admissible isomorphism $T^H \rightarrow T$ defined over $F$ and mapping $t^H$ to $t$. This definition is the same as in
\cite{LS90}, but our wording is opposite -- in \cite{LS90} the element $t^H$ is called an image of $t$. If $t$ is an image of $t^H$ we will also call $(t^H,t)$
a pair of related elements. For such a pair, we consider the set of $\phi : T^H \rightarrow T$, where $T^H$ is a maximal torus in $H$ containing $t^H$, $T$ is
a maximal torus of $G$ containing $t$, and $\phi$ is an admissible isomorphism defined over $F$ and mapping $t^H$ to $t$. On this set we define an equivalence
relation, by saying that two such isomorphisms $\phi$ and $\phi'$ are $(G_t,H_{t^H})$-equivalent if there exist $g \in G_t(\ol{F})$ and $h \in H_{t^H}(\ol{F})$
s.t. $\phi' = \tx{Ad}(g)\phi\tx{Ad}(h)$. If $\phi$ is an element of this set, and $H_{t^H}$ is quasi-split, then $H_{t^H}$ can be identified with an endoscopic
group of $G_t$ in such a way that $\phi$ becomes an admissible isomorphism with respect to $(G_t,H_{t^H})$. Then we can talk about images, admissible
isomorphisms, etc. with respect to the group $G_t$ and its endoscopic group $H_{t^H}$. When we do so, we will use the prefix $(G_t,H_{t^H},\phi)$.

If $\gamma,\gamma'$ are two strongly $G$-regular semi-simple elements, each of which belongs to either $G(F)$ or $H(F)$, and $T,T'$ are their centralizers,
then there exists at most one admissible isomorphism $T \rightarrow T'$ which maps $\gamma$ to $\gamma'$. We will call this isomorphism
$\phi_{\gamma,\gamma'}$. If it exists, then so does $\phi_{\gamma',\gamma}$ and $\phi_{\gamma',\gamma}=\phi_{\gamma,\gamma'}^{-1}$. Moreover, if
$\gamma,\gamma',\gamma''$ are three elements as above and $\phi_{\gamma,\gamma'}$ and $\phi_{\gamma',\gamma''}$ exist, then so does $\phi_{\gamma,\gamma''}$
and
\[ \phi_{\gamma,\gamma''} = \phi_{\gamma',\gamma''} \circ \phi_{\gamma,\gamma'} \]
The same can also be done with regular semi-simple elements of the Lie algebras of $G$ and $H$ and we will use the same notation for that case.

\section{Pure inner twists} \label{sec:pure}

Let $A,B$ be reductive groups over $F$. A pure inner twist
\[ (\psi,z) : A \rightarrow B \]
consists of an isomorphism of $\ol{F}$-groups $\psi : A \times \ol{F} \rightarrow B \times \ol{F}$ and an element $z \in Z^1(\Gamma,A)$ s.t.
\[ \forall \sigma \in \Gamma : \psi^{-1}\sigma(\psi) = \tx{Ad}(z_\sigma) \]
We will from now on abbreviate "pure inner twist" to simply "twist", since these will be the only twists of reductive groups that will concern us here.

The twist $(\psi,z)$ is called trivial the image of $z$ in $H^1(\Gamma,A)$ is trivial. In that case there exists $a \in A(\ol{F})$ s.t.
\[ \psi\circ\tx{Ad}(a) : A \rightarrow B \]
is an isomorphism over $F$. Clearly the element $a$ is unique up to right multiplication by $A(F)$. We will call the twist $(\psi,z)$ strongly trivial if
$z=1$. In that case of course $\psi$ is already defined over $F$. An example of a trivial twist is given by $(\tx{Ad}(g),g^{-1}\sigma(g)) : A \rightarrow A$
for any $g \in A(\ol{F})$. This twist is strongly trivial if and only if $g \in A(F)$.

Starting from $(\psi,z) : A \rightarrow B$ we can form the inverse twist $(\psi,z)^{-1} : B \rightarrow A$, which is given by
$(\psi^{-1},\psi(z_\sigma^{-1}))$.

If $(\psi,z) : A \rightarrow B$ and $(\phi,u) : B \rightarrow C$ are twists, then we can form their composition
\[(\phi,u) \circ (\psi,z) : A \rightarrow C \]
which is given by $(\phi\circ\psi,\psi^{-1}(u)z)$. One immediately checks
\begin{eqnarray*}
(\psi,z) \circ (\psi,z)^{-1}&=&(\tx{id}_B,1)\\
(\psi,z)^{-1} \circ (\psi,z)&=&(\tx{id}_A,1)\\
{[(\phi,u) \circ (\psi,z)]}^{-1}&=&(\psi,z)^{-1} \circ (\phi,u)^{-1}\\
(\chi,v)\circ [(\phi,u)\circ(\psi,z)] &=&[(\chi,v)\circ(\phi,u)]\circ(\psi,z)
\end{eqnarray*}
In particular, reductive groups and pure inner twists form a groupoid.

Let $(\psi,z), (\psi',z') : A \rightarrow B$ be two twists. They are called equivalent if $(\psi',z') \circ (\psi,z)^{-1}$ equals
$(\tx{Ad}(g),g^{-1}\sigma(g))$ for some $g \in B(\ol{F})$. One immediately checks the equality
\[ (\psi,z)^{-1} \circ (\tx{Ad}(g),g^{-1}\sigma(g)) \circ (\psi,z) = (\tx{Ad}(h),h^{-1}\sigma(h)),\quad h=\psi^{-1}(g) \]
from which it follows that this defines an equivalence relation on all inner twists which is invariant under composition and taking inverses.


\subsection{Conjugacy along pure inner twists}

Now consider a twist $(\psi,z) : A \rightarrow B$ and two elements $a \in A(F), b \in B(F)$. We call $a,b$ conjugate (with respect to $(\psi,z)$) if there
exists a twist $(\psi',z')$ equivalent to $(\psi,z)$ which maps $a$ to $b$ and is strongly trivial. We call $a,b$ stably conjugate (with respect to $(\psi,z)$)
if there exists a twist $(\psi',z')$ equivalent to $(\psi,z)$ which maps $a$ to $b$ and descends to a twist $A_a \rightarrow B_b$. The latter condition simply
means that $z'$ takes values in $B_b$ (a-priori it only takes values in $\tx{Cent}(b,B)$).

The following is immediately clear
\begin{fct}\ \\[-20pt]
\begin{enumerate}
\item Applied to the twist $(\tx{id},1) : A \rightarrow A$ the notions defined above coincide with the usual ones for the group A
\item If $a \in A(F),b\in B(F)$ are conjugate with respect to $(\psi,z) : A \rightarrow B$, then they are also stably conjugate and moreover $(\psi,z)$ is
    a trivial twist
\item If $a\in A(F),b\in B(F)$ are conjugate (resp. stably conjugate) with respect to $(\psi,z) : A \rightarrow B$, then so are they with respect to any
    twist equivalent to $(\psi,z)$.
\item If $a \in A(F)$ and $b \in B(F)$ are conjugate (resp. stably-conjugate) with respect to $(\psi,z) : A \rightarrow B$, then so are they with respect
    to $(\psi,z)^{-1} : B \rightarrow A$
\item If $(\psi,z): A \rightarrow B$ and $(\phi,u) : B \rightarrow C$ are two twists and $a \in A(F),b\in B(F),c\in C(F)$ are s.t. $a,b$ and $b,c$ are
    conjugate (resp. stably-conjugate), then so are $a,c$.
\end{enumerate}
\end{fct}

Let $a \in A(F)$ and $b \in B(F)$ be stably conjugate assume that $\tx{Cent}(a,A)$ is connected. Choose a twist $(\phi,u) : A \rightarrow B$ which is
equivalent to $(\psi,z)$ and sends $a$ to $b$, and write $\tx{inv}(a,b)$ for the image of $u$ in $H^1(\Gamma,A_a)$.

\begin{fct}\ \\[-20pt]
\begin{enumerate}
\item The element $\tx{inv}(a,b)$ is independent of the choice of the twist $(\phi,u)$.
\item Applied to the twist $(\tx{id},1) : A \rightarrow A$, $\tx{inv}$ coincides with the usual definition for the group $A$.
\item The image of $\tx{inv}(a,b)$ in $H^1(\Gamma,A)$ equals $z$.
\end{enumerate}
\end{fct}
\pf This is obvious.\qed

\begin{fct} \label{fct:invprod} Let $a \in A(F)$, $b \in B(F)$ and $c \in C(F)$ be s.t. the inner twists
\[ A \stackrel{(\phi,u)}{\longrightarrow} B \stackrel{(\psi,z)}{\longrightarrow} C \]
send $a$ to $b$ to $c$. Assume that $\tx{Cent}(a,A)$ is connected. Then
\[ \tx{inv}(a,c) = \phi^{-1}(\tx{inv}(b,c))\tx{inv}(a,b) \]
\end{fct}
\pf This follows at once from the composition formula for twists.\qed

Now let $A$ be quasi-split. We consider a set $I$ of triples $(A^z,\psi_z,z)$ s.t. $(\psi_z,z) : A \rightarrow A^z$ is a twist. Put
\[ A^I = \bigsqcup_{(A^z,\psi_z,z)} A^z \]
This is a variety over $F$ (it will not be of finite type if $I$ is infinite). For $a \in A^z$ and $b \in A^{z'}$ we obtain notions of conjugacy and stable
conjugacy, namely those relative to the twist $(\psi_{z'},z') \circ (\psi_z,z)^{-1}$. Thus we can talk about conjugacy classes and stable conjugacy classes of
elements of $A^I(F)$. For the sake of abbreviation, we will call a twist $(\phi,u) : A^z \rightarrow A^{z'}$ allowable if it is equivalent to $(\psi_{z'},z')
\circ (\psi_z,z)^{-1}$. Note that the set of allowable twists is invariant under composing and taking inverses.

\begin{fct} \label{fct:stcqs} Every stable conjugacy class of $A^I(F)$ meets $A(F)$. \end{fct}
\pf This is a consequence of a well known theorem of Steinberg, which implies that any maximal torus of a reductive group transfers to its quasi-split inner
form.\qed

\begin{lem} \label{lem:stcit}
Let $\bar I$ be the image of $I$ in $H^1(\Gamma,A)$ under the map $(A^z,\psi_z,z) \mapsto [z]$. Then for each $a \in A(F)$ whose centralizer is connected, the
map $b \mapsto \tx{inv}(a,b)$ is a bijection from the set of conjugacy classes inside the stable class of $a$ in $A^I(F)$ to the preimage of
    $\bar I$ under $H^1(\Gamma,A_a) \rightarrow H^1(\Gamma,A)$.
\end{lem}

\rmk One can prove a similar lemma for $a \in A^{z'}(F)$ and any $z'$, but the statement is more awkward and we will not need it.

\pf Let $b \in A^z(F)$ and $b' \in A^{z'}(F)$ be conjugate elements belonging to the stable class of $a$. Thus there exists an allowable strongly trivial twist
$(\chi,1) : A^{z} \rightarrow A^{z'}$ mapping $b$ to $b'$. Let $(\phi,u) : A \rightarrow A^z$ be an allowable twist mapping $a$ to $b$, thus
$\tx{inv}(a,b)=[u]$. Then $(\chi,1) \circ (\phi,u)$ is an allowable twist $A \rightarrow A^{z'}$, mapping $a$ to $b'$, so $\tx{inv}(a,b')$ equals the class of
the cocycle of $[(\chi,1)\circ(\phi,u)]$, which is also $[u]$. This shows that $\tx{inv}(a,b)=\tx{inv}(a,b')$ and we see that the map $b \mapsto \tx{inv}(a,b)$
is a well-defined map on the set of conjugacy classes inside the stable class of $a$. By above facts it lands in the preimage of $\bar I$. We will show that it
is injective. To that end, let $b \in A^z(F)$ and $b' \in A^{z'}(F)$ be s.t. $\tx{inv}(a,b)=\tx{inv}(a,b')$. Let $(\phi,u) : A \rightarrow A^z$ and $(\phi',u')
: A \rightarrow A^{z'}$ be allowable twists sending $a$ to $b$ reps. $b'$. By assumption there exists $i \in A_a$ s.t. $u=i^{-1}u'\sigma(i)$. But then
$(\phi',u') \circ (\tx{Ad}(i),i^{-1}\sigma(i))$ is again an allowable twist $A \rightarrow A^{z'}$ sending $a$ to $b'$, and so replacing $(\phi',u')$ by it we
achieve $u=u'$. But now it is clear that $(\phi',u')\circ(\phi,u)^{-1}$ is an allowable strongly trivial twist $A^z \rightarrow A^{z'}$ sending $b$ to $b'$,
thus showing that $b$ and $b'$ are conjugate. Finally we show that the map $b \mapsto \tx{inv}(a,b)$ is surjective. Thus let $[u] \in H^1(\Gamma,A_a)$ be an
element mapping to $[z] \in H^1(\Gamma,A)$, where $(A^z,\psi_z,z) \in I$. Then there exists $g \in A$ s.t. $u=g^{-1}z\sigma(g)$. Put $b = \psi_z(\tx{Ad}(g)a)$.
One computes immediately that $b \in A^z(F)$. By construction $(\psi_z,z) \circ (\tx{Ad}(g),g^{-1}\sigma(g))$ maps $a$ to $b$, which shows that $a$ and $b$ are
stably conjugate and that $\tx{inv}(a,b)=g^{-1}z\sigma(g)$.\qed


\subsection{Transfer factors for pure inner twists} \label{sec:pure_tf}

Let $G$ be a quasi-split $F$-group, $(\psi,z) : G \rightarrow G'$ a twist, and $(H,s,{^L\eta})$ an extended triple for $G$. Then $(H,s,{^L\eta})$ is also an
extended triple for $G'$. This data gives canonical relative geometric transfer factors $\Delta^G_H(\gamma^H,\gamma,\bar\gamma^H,\bar\gamma)$ for $(G,H)$ and
$\Delta^{G'}_H(\gamma^H,\gamma',\bar\gamma^H,\bar\gamma')$ for $(G',H)$ (see \cite{LS87}). Let $\Delta^G_H(\gamma^H,\gamma)$ be an arbitrary normalization for
the absolute transfer factor for $(G,H)$. For any pair $\gamma^H \in H(F)$ and $\gamma' \in G'(F)$ of strongly $G$-regular related elements we choose an
element $\gamma \in G(F)$ stably conjugate to $\gamma'$ (which exists by Fact \ref{fct:stcqs}) and define
\[ \Delta^{G'}_H(\gamma^H,\gamma') = \Delta^G_H(\gamma^H,\gamma) \cdot \langle\tx{inv}(\gamma,\gamma'),\hat\phi_{\gamma,\gamma^H}(s)\rangle^{-1} \]
where
\[ \langle\rangle : H^1(\Gamma,T) \times \pi_0(\hat T^\Gamma) \rightarrow \C^\times \]
is the Tate-Nakayama pairing, and $T=\tx{Cent}(\gamma,G)$.

\begin{lem} $\Delta^{G'}_H(\cdot,\cdot)$ is well defined and is an absolute transfer factor for $(G',H)$ \end{lem}
\pf We need to show that $\Delta^{G'}_H(\gamma^H,\gamma')$ is independent of the choice of $\gamma$. Thus let $\tilde\gamma \in G(F)$ be another element in the
stable class of $\gamma'$. We know from \cite{LS87}
\[ \Delta^G_H(\gamma^H,\tilde\gamma)=\Delta^G_H(\gamma^H,\gamma)\langle\tx{inv}(\gamma,\tilde\gamma),\hat\phi_{\gamma,\gamma^H}(s)\rangle^{-1} \]
On the other hand if $(\phi,u) : A \rightarrow A$ is an admissible twist mapping $\tilde\gamma$ to $\gamma$, then
$\phi_{\tilde\gamma,\gamma^H}=\phi_{\gamma,\gamma^H}\circ\phi$ and by functoriality of the Tate-Nakayama pairing we get
\[ \langle\tx{inv}(\tilde\gamma,\gamma'),\hat\phi_{\tilde\gamma,\gamma^H}(s)\rangle^{-1} =
\langle\phi(\tx{inv}(\tilde\gamma,\gamma')),\hat\phi_{\gamma,\gamma^H}(s)\rangle^{-1} \]

Thus
\begin{eqnarray*}
&&\Delta^G_H(\gamma^H,\tilde\gamma)\langle\tx{inv}(\tilde\gamma,\gamma'),\hat\phi_{\tilde\gamma,\gamma^H}(s)\rangle^{-1}=\\
&&\Delta^G_H(\gamma^H,\gamma)\langle\tx{inv}(\gamma,\tilde\gamma),\hat\phi_{\gamma,\gamma^H}(s)\rangle^{-1}\
\langle\phi(\tx{inv}(\tilde\gamma,\gamma')),\hat\phi_{\gamma,\gamma^H}(s)\rangle^{-1}=\\
&&\Delta^G_H(\gamma^H,\gamma)\langle\tx{inv}(\gamma,\gamma'),\hat\phi_{\gamma,\gamma^H}(s)\rangle^{-1}
\end{eqnarray*}
by Fact \ref{fct:invprod}.

This shows that $\Delta^{G'}_H(\gamma^H,\gamma')$ is independent of the choice of $\gamma$. To show that it is an absolute transfer factor for $(G',H)$ we must
prove for any two strongly $G$-regular related pairs $(\gamma^H,\gamma')$ and $(\bar\gamma^H,\bar\gamma')$ in $H(F)\times G'(F)$ the equality
\[ \frac{\Delta^{G'}_H(\gamma^H,\gamma')}{\Delta^{G'}_H(\bar\gamma^H,\bar\gamma')} = \Delta^{G'}_H(\gamma^H,\gamma',\bar\gamma^H,\bar\gamma') \]
which by construction of $\Delta^{G'}_H$ is equivalent to
\[ \frac{\langle\tx{inv}(\gamma,\gamma'),\hat\phi_{\gamma,\gamma^H}(s)\rangle^{-1}}{\langle\tx{inv}(\bar\gamma,\bar\gamma'),\hat\phi_{\bar\gamma,\bar\gamma^H}(s)\rangle^{-1}}=
   \frac{\Delta^{G'}_H(\gamma^H,\gamma',\bar\gamma^H,\bar\gamma')}{\Delta^{G}_H(\gamma^H,\gamma,\bar\gamma^H,\bar\gamma)}    \]
where $\gamma \in G(F)$ is any element in the stable class of $\gamma'$, and $\bar\gamma\in G(F)$ is any element in the stable class of $\bar\gamma'$. Applying
\cite[Lemma 4.2.A]{LS87} we need to show
\[ \frac{\langle\tx{inv}(\gamma,\gamma'),\hat\phi_{\gamma,\gamma^H}(s)\rangle^{-1}}{\langle\tx{inv}(\bar\gamma,\bar\gamma'),\hat\phi_{\bar\gamma,\bar\gamma^H}(s)\rangle^{-1}}=
   \left\langle \tx{inv}\left(\frac{\gamma,\gamma'}{\bar\gamma,\bar\gamma'}\right),s_U \right\rangle \]
The right hand side of this equality is constructed in \cite[\S3.4]{LS87}. Working through the construction, one sees that in our case the objects are as
follows: Let $T$ and $\ol{T}$ denote the centralizers of $\gamma$ and $\bar\gamma$ in $G$, and let $T_\tx{sc}$ and $\ol{T}_\tx{sc}$ be their preimages in the
simply connected cover $G_\tx{sc}$ of the derived group of $G$. Let $Z_\tx{sc}$ be the center of $G_\tx{sc}$. Then $U=T_\tx{sc} \times \ol{T}_\tx{sc} /
\{(z^{-1},z)|\ z \in Z_\tx{sc} \}$. We have the following dual diagrams
\begin{diagram}
T \times \ol{T}&\lTo&T_\tx{sc} \times \ol{T}_\tx{sc}&\rTo&U\\
\hat T \times \hat{\ol{T}}&\rTo&\hat T_\tx{ad} \times \hat{\ol{T}}_\tx{ad}&\lTo&\hat U
\end{diagram}

The elements $s_U \in \hat U^\Gamma$ and $(\hat\phi_{\gamma^H,\gamma}(s),\hat\phi_{\bar\gamma^H,\bar\gamma}(s)) \in \hat T \times \hat{\ol{T}}$ map to the same
element in $\hat T_\tx{ad} \times \hat{\ol{T}}_\tx{ad}$. There exists an element of $H^1(\Gamma,T_\tx{sc}) \times H^1(\Gamma,\ol{T}_\tx{sc})$ which maps to
$\tx{inv}\left(\frac{\gamma,\gamma'}{\bar\gamma,\bar\gamma'}\right) \in H^1(\Gamma,U)$ and to $(\tx{inv}(\gamma,\gamma')^{-1},\tx{inv}(\bar\gamma,\bar\gamma'))
\in H^1(\Gamma,T) \times H^1(\Gamma,\ol{T})$. The equality now follows again from the functoriality of the Tate-Nakayama pairing.\qed

Now let $I$ be a set of pure inner twists for $G$ and construct $G^I$ as above. Taking the disjoint union over $I$ of all functions $\Delta^{G^z}_H$ we obtain
a function
\[ \Delta^{G^I}_H : H_{G-\tx{sr}}(F) \times G^I_\tx{sr}(F) \longrightarrow \C^\times \]
\begin{fct} For all stably conjugate $\gamma,\gamma' \in G^I_\tx{sr}(F)$ and $\gamma^H \in H_{G-\tx{sr}}(F)$ we have
\[ \Delta^{G^I}_H(\gamma^H,\gamma') = \Delta^{G^I}_H(\gamma^H,\gamma) \cdot \langle \tx{inv}(\gamma,\gamma'),\hat\phi_{\gamma,\gamma^H}(s) \rangle^{-1} \]
\end{fct}
\pf Let $\gamma_0 \in G(F)$ be an element stably conjugate to $\gamma$ (it exists by Fact \ref{fct:stcqs}). Then by construction of $\Delta^{G^I}_H$ we have
\[ \Delta^{G^I}_H(\gamma^H,\gamma')\Delta^{G^I}_H(\gamma^H,\gamma)^{-1} =
\langle\tx{inv}(\gamma_0,\gamma')\tx{inv}(\gamma_0,\gamma)^{-1},\hat\phi_{\gamma_0,\gamma_H}(s)\rangle^{-1} \]

By Fact \ref{fct:invprod} the right hand side equals $\langle \phi_{\gamma,\gamma_0}(\tx{inv}(\gamma,\gamma')),\hat\phi_{\gamma_0,\gamma_H}(s)\rangle^{-1}$ and
the claim now follows from the functoriality of the Tate-Nakayama pairing.\qed

\rmk We see in particular the the function $\gamma \mapsto \Delta^{G^I}_H(\gamma^H,\gamma)$ is constant on the conjugacy classes of $G^I(F)$.


\section{Statement of the main result} \label{sec:main}

We fix an unramified reductive group $G$ over $F$, and a Borel pair $(T_0,B_0)$ of $G$ defined over $F$. Then $\Gamma$ acts on $X^*(T_0)$ through a finite
cyclic subgroup of $\tx{Aut}(X^*(T_0))$ generated by the image of $\tx{Fi}$; we will denote by $\vartheta$ both this image as well as its dual in
$\tx{Aut}(X_*(T_0))$. Let $(\hat G,\hat B_0,\hat T_0)$ be the dual datum to $(G,B_0,T_0)$. If $\Omega(T_0,G)$ and $\Omega(\hat T_0,\hat G)$ denote the
corresponding Weyl-groups, then there is a natural isomorphism between them given by duality. We choose an $L$-group $^LG$ for $G$ s.t. the $\Gamma$-action on
$\hat G$ preserves the pair $(\hat B_0,\hat T_0)$.

We also fix an endoscopic triple $(H,s,\hat\eta)$ for $G$ s.t. $H$ is unramified. We choose again a Borel pair $(T_0^H,B_0^H)$ defined over $F$, let $(\hat
H,\hat B_0^H,\hat T_0^H)$ be the dual datum to $(H,B_0^H,T_0^H)$ and ${^LH}$ an $L$-group for $H$ s.t. the $\Gamma$-action on $\hat H$ preserves $(\hat
B_0^H,\hat T_0^H)$.

We choose a hyperspecial point $o$ in the apartment of $T_0$ and obtain an $O_F$-structure on $G$ and $\mf{g}$. Then $G_o,G_{o^+}$ resp.
$\mf{g}_o,\mf{g}_{o^+}$ will be the parahoric and its pro-unipotent radical of $G(O_{F^u})$ resp. $\mf{g}(O_{F^u})$ associated to $o$. We also choose a
hyperspecial point, denoted again by $o$, in the apartment of $T_0^H$ and obtain the same structures on $H$ and $\mf{h}$.

Up to equivalence the map $\hat\eta : \hat H \rightarrow \hat G$ may be chosen so that $\hat\eta^{-1}(\hat B_0,\hat T_0) = (\hat B_0^H,\hat T_0^H)$. Then we
have in particular an isomorphism of complex tori $\hat\eta_{T_0^H} : \hat T_0^H \rightarrow \hat T_0$. There exists an element $\omega \in
Z^1(\Gamma,\Omega(\hat T_0,\hat G))$ s.t. $\omega(\sigma)\sigma\circ\hat\eta|_{\hat T_0^H}\circ\sigma^{-1}=\hat\eta|_{\hat T_0^H}$ for all $\sigma \in \Gamma$.
Thus we dually obtain an isomorphism of $F$-tori $\eta : T_0^\omega \rightarrow T_0^H$, where $T_0^\omega$ denotes the twist of $T_0$ by $\omega$.

By \cite[Lemma 6.1]{Hal93} the map $\hat\eta : \hat H \rightarrow \hat G$ can be extended to an $L$-embedding ${^L\eta} : {^LH} \rightarrow {^LG}$ in such a
way, that the 1-cocycle
\[ I_F \rightarrow {^LH} \rightarrow \hat H \]
is trivial. We choose such an extension. The extended triple $(H,s,{^L\eta})$ is then unramified in the sense of \cite{Hal93}.


\subsection{Review of the construction of DeBacker and Reeder} \label{sec:review}

In this section we want to review the construction from \cite{DR09} of the $L$-packet on $G$ and its pure inner forms corresponding to a Langlands parameter
$\phi : W_F \rightarrow {^LG}$ which is TRSELP in the sense of loc. cit. Our purpose is not to review the details of the construction, but rather to gather the
necessary notation and properties needed in the subsequent sections.

Recall that $\phi$ is called TRSELP if it is trivial on $\tx{SL}_2(\C)$, $\tx{Cent}(\phi(I_F),\hat G)$ is a maximal torus of $\hat G$, and $Z(\hat G)^\Gamma$
is of finite index in $\tx{Cent}(\phi,\hat G)$. Up to equivalence we may assume that $\phi(I_F) \subset \hat T_0$. There is an element $w \in
Z^1(\Gamma,\Omega(\hat T_0,\hat G))$ s.t.
\begin{eqnarray*}
\tx{Ad}(\phi(\sigma))|_{\hat T_0}&=&w(\sigma)\sigma,\quad \forall \sigma \in W_F
\end{eqnarray*}

Let $T_0^w$ be the twist of $T_0$ by $w$. The ellipticity of $\phi$ implies that $T_0^w/Z$ is anisotropic. Put $X=X_*(T_0^w)$. This is a $\Z[\Gamma]$-module,
where the $\Gamma$-action comes from that on $T_0^w$. Let $\bar X$ be the quotient of $X$ by the coroot-lattice, and $X_\Gamma$ resp. $\bar X_\Gamma$ denote
the $\Gamma$-coinvaraints in $X$ resp. $\bar X$. Let $X_w$ be the preimage of $[X_\Gamma]_\tx{tor}$ in $X$. Write $C_\phi$ for the component group of the
centralizer of $\phi$ in $\hat G$. We have the following diagram
\begin{diagram}[LaTeXeqno]\label{diag}
\tx{Irr}(C_\phi)&\rEquals&\tx{Irr}(\pi_0(\hat{T_0^w}^\Gamma))&\rTo&\tx{Irr}(\pi_0(Z(\hat G)^\Gamma)) \\
&&\uTo<{\cong}&&\uTo>{\cong}\\
X_w&\rTo&[X_\Gamma]_\tx{tor}&\rTo&[\bar X_\Gamma]_\tx{tor}\\
&&\dTo<{\cong}>{\tx{DR}_{T_0^w}}&&\dTo>{\cong}<{\tx{DR}_G}\\
&&H^1(\Gamma,T_0^w)&\rTo&H^1(\Gamma,G)
\end{diagram}
The bottom square of it is \cite[Lemma 2.6.1]{DR09}. The top equality follows from
\[\tx{Cent}(\phi,\hat G) = \hat T_0^{w\vartheta} = \hat{T_0^w}^\Gamma\]
while the rest is given by the obvious restriction maps.

The map $X_w \rightarrow H^1(\Gamma,G)$ in this diagram will be denoted by $r$. For any $u \in H^1(\Gamma,G)$ we let $[r^{-1}(u)]$ be the image of $r^{-1}(u)$
in $[X_\Gamma]_\tx{tor}$. The map $X_w \rightarrow \tx{Irr}(C_\phi)$ will be denoted by $\lambda \mapsto \rho_\lambda$.

From the Langlands parameter $\phi$ DeBacker and Reeder construct (see \cite[\S4]{DR09}) a Langlands parameter $\phi_T : W_F \rightarrow {^LT_0^w}$ which
corresponds to a regular depth-zero character $\theta : T_0^w(F) \rightarrow \C^\times$ (both notations $\theta$ and $\chi_\phi$ are used for this character in
loc.cit). Moreover, given $\lambda \in X_w$, they construct the following objects

\begin{itemize}
\item An element $u_\lambda \in Z^1(\Gamma,G)$ (trivial on inertia). Let $(\psi_\lambda,u_\lambda) : G \rightarrow G^\lambda$ be the corresponding twist.
\item A maximal torus $T_\lambda \subset G^\lambda$, together with an element $p_\lambda \in G^\lambda(F^u)$ s.t.
\[ \tx{Ad}(p_\lambda)\psi_\lambda : T_0^w \rightarrow T_\lambda \]
is an isomorphism of $F$-tori
\item A depth-zero supercuspidal representation $\pi_\lambda$ of $G^\lambda(F)$.
\end{itemize}

Furthermore they show in the proof of \cite[Thm 4.5.3]{DR09} that for $\lambda,\mu \in X_w$ one has $\rho_\lambda=\rho_\mu$ if and only if $\psi_\lambda \circ
\psi_\mu^{-1} : G^\mu \rightarrow G^\lambda$ is a trivial twist and the transfer of $\pi_\mu$ to $G^\lambda$ with respect to one (hence any) strongly trivial
twist equivalent to $\psi_\lambda \circ \psi_\mu^{-1}$ coincides with $\pi_\lambda$. Thus if we put
\[I = \{ (G^\lambda,\psi_\lambda,u_\lambda)|\ \lambda \in X_w \} \]
and construct $G^I$ as in Section \ref{sec:pure}, then for each $\rho \in \tx{Irr}(C_\phi)$ we obtain a conjugation\-invariant function $\Theta_{\phi,\rho}$ on
$G^I(F)$ by taking any $\lambda \in X_w$ s.t. $\rho_\lambda=\rho$ and extending the character of $\pi_\lambda$ to a conjugation\-invariant function on
$G^I(F)$.

To simplify their stability calculations, DeBacker and Reeder rigidify their constructions in the following way. In every class of $H^1(\Gamma,G)$ they choose
a specific representative $u \in Z^1(\Gamma,G)$, which again gives rise to a twist $(\psi,u) : G \rightarrow G^u$. For each $\lambda \in r^{-1}(u)$ they
construct an element $q_\lambda \in G^u(F^u)$ s.t. the maximal torus $S_\lambda = \tx{Ad}(q_\lambda)\psi(T_0)$ is defined over $F$ and
\[ \tx{Ad}(q_\lambda)\psi : T_0^w \rightarrow S_\lambda \]
is an isomorphism over $F$. For any strongly regular semi-simple element $Q \in S_0(F)$ the map
\[ \lambda \mapsto \tx{Ad}(q_\lambda)\psi\tx{Ad}(q_0^{-1})Q \]
is a bijection from $[r^{-1}(u)]$ to a set of representatives for the stable class of $Q$ in $G^u(F)$ (\cite[Lem. 2.10.1]{DR09}). In particular, the tori
$S_\lambda$ exhaust the stable class of $T_0^w$ in $G^u$. It will be important for later to note that $p_0=q_0 \in G(O_{F^u})$. For every $\rho \in
\tx{Irr}(C_\phi)$ mapping to the class of $u$, they define a representation $\pi_u(\phi,\rho)$ on $G^u(F)$. It is equal to the transfer of $\pi_\lambda$ via
any strongly trivial twist $G^\lambda \rightarrow G^u$ equivalent to $\psi \circ \psi_\lambda^{-1}$, where $\lambda$ is any element of $r^{-1}(u)$.

It is clear from the constructions that for any $\lambda \in r^{-1}(u)$, the twist $\psi_\lambda \circ \psi^{-1}$ defines an injection from the conjugacy
classes in $G^u(F)$ to the conjugacy classes in $G^I(F)$ whose image consists of those conjugacy classes which meet $G^\mu(F)$ for $\mu \in r^{-1}(u)$.
Moreover, this twist identifies the character of $\pi_u(\phi,\rho)$ with the function $\Theta_{\phi,\rho}$, where both are viewed as class functions.

The same construction can be applied to a TRSELP $\phi^H : W_F \rightarrow {^LH}$ and the corresponding objects will carry the superscript $H$.


\subsection{The Whittaker character}

We extend the chosen pair $(T_0,B_0)$ of $G$ to a splitting $(T_0,B_0,\{X_\alpha\})$ where each simple root vector $X_\alpha$ is chosen so that the
homomorphism
\[ \mb{G}_a \rightarrow G \]
determined by it is defined over $O_{F^u}$ and the image of $1$ under
\[ \mb{G}_a(O_{F^u}) \rightarrow G(O_{F^u}) \rightarrow G(\ol{\mb{F}_q}) \]
is non-trivial. Such a splitting is called admissible by \cite{Hal93}. Let $N$ denote the unipotent radical of $B_0$.

\begin{lem} There exists an additive character $\psi : F \rightarrow \C^\times$ which is non-trivial on $O_F$ but trivial on $\pi O_F$, s.t.
the representation $\pi_1(\phi,1)$ is generic with respect to the character $N(F) \rightarrow \C^\times$ determined by $\psi$ and the chosen splitting.
\end{lem}

\pf The representation $\pi_1(\phi,1)$ is the same as the representation $\pi_0$ defined in \cite[\S4.5]{DR09}. By Lemmas 6.2.1 and 6.1.2 in loc. cit. it is
generic with respect to a character $N(F) \rightarrow \C^\times$ which has depth-zero at $o$. This character is generic and is thus given by the composition of
the $F$-homomorphism
\[ N \rightarrow \prod_{\alpha \in \Delta} \mb{G}_a \stackrel{\Sigma}{\longrightarrow} \mb{G}_a \]
determined by the chosen splitting with an additive character
\[ \psi : F \rightarrow \C^\times \]
The choice of the simple root vectors $X_\alpha$ ensures that the homomorphism $N \rightarrow \mb{G}_a$ is in fact defined over $O_F$ and moreover maps
$N(O_F)$ surjectively onto $\mb{G}_a(O_F)$. The genericity of the character $N(F) \rightarrow \C^\times$ now implies that $\psi$ is non-trivial on $O_F$ and
trivial on $\pi O_F$. \qed

From now on we fix an additive character $\psi : F \rightarrow \C^\times$ as in the above Lemma.


\subsection{Definition of the unstable character} \label{sec:schar_def}

For $t \in \tx{Cent}(\phi,\hat G)$ we define on $G^I(F)$ the function
\[ \Theta_\phi^{t} = \sum_{\rho \in \tx{Irr}(C_\phi)} e_\rho\tr\rho(t)\Theta_{\phi,\rho} \]
where for any $\lambda \in X_w$ with $\rho_\lambda=\rho$ we put $e_\rho=e(G^\lambda)$, the latter being the sign defined in \cite{Kot83}. This is the
$t$-unstable character corresponding to the packet $\Pi(\phi)$ defined in \cite[\S4.5]{DR09}.

We will also define the $t$-unstable character of the normalized $L$-packet $\Pi_u(\phi)$ defined in \cite[\S4.6]{DR09} for the specific twists $(\psi,u) : G
\rightarrow G^u$ considered there. This character is
\[ \Theta^{t}_{\phi,u} := e(G^u)\sum_{\rho \in \tx{Irr}(C_\phi,u)} \tr\rho(t) \Theta_{\pi_u(\phi,\rho)} \]
where $\tx{Irr}(C_\phi,u)$ is the fiber over $u$ of the map $\tx{Irr}(C_\phi) \rightarrow H^1(\Gamma,G)$ given in diagram \eqref{diag}. We will show in Lemma
\ref{lem:kotdr} that the map $H^1(\Gamma,G) \rightarrow \pi_0(Z(\hat G)^\Gamma)$ in diagram \eqref{diag} is a particular normalization of the Kottwitz
isomorphism, and so the set $\tx{Irr}(C_\phi,u)$ is the set of all irreducible representations of $C_\phi$ which transform under $\pi_0(Z(\hat G)^\Gamma)$ by
the character corresponding to $u$ via the Kottwitz isomorphism.

The restriction of $\Theta^1_\phi$ to $G(F)$, which also equals $\Theta^1_{\phi,1}$, will be denoted by $\mc{S}\Theta_\phi$.


\subsection{Statement of the main result} \label{sec:mainst}

Before stating the main result, we need to impose some mild conditions on the residual characteristic of $F$. These restrictions come from the papers
\cite{DR09} and \cite{Hal93}. To state them, let $n_G$ denote the smallest dimension of a faithful representation of $G$, and $n_H$ be the corresponding number
for $H$. Let $e$ be the ramification degree of $F/\Q_p$ and $e_G$ be the minimum over the ramification degrees (again over $\Q_p$) of all splitting fields of
maximal tori of $G$. The restrictions we impose are

\begin{itemize}
\item $q_F \geq |R(T_0,B_0)|$
\item $p \geq (2+e)\max(n_G,n_H)$
\item $p\geq 2+e_G$
\end{itemize}
The first two items are imposed in \cite[\S12.4]{DR09}, while the third is imposed in the main result of \cite{Hal93} -- Theorem 10.18.

From now on we assume that these restrictions hold.

Let $\phi^H : W_F \rightarrow {^LH}$ be a Langlands parameter for $H$, then $\phi = {^L\eta} \circ \phi^H$ is a Langlands parameter for $G$. We are interested
in the situation in which both $\phi^H$ and $\phi$ are TRSELP. Then $(H,s,\hat\eta)$ is automatically an elliptic endoscopic triple for $G$. Up to equivalence
we may assume that $\phi^H$ maps inertia into $\hat T_0^H$, then $\phi$ maps inertia into $\hat T_0$ by our choice of $\hat\eta$. There are elements $w \in
Z^1(\Gamma,\Omega(\hat T_0,\hat G))$, $w^H \in Z^1(\Gamma,\Omega(\hat T_0^H,\hat H))$ s.t.
\begin{eqnarray*}
\tx{Ad}(\phi(\sigma))|_{\hat T_0}&=&w(\sigma)\sigma,\quad \forall \sigma \in W_F\\
\tx{Ad}(\phi^H(\sigma))|_{\hat T_0^H}&=&w^H(\sigma)\sigma,\quad \forall \sigma \in W_F
\end{eqnarray*}

Let $\Delta_{\psi}$ be the Whittaker normalization \cite[\S5.3]{KS99} of the absolute transfer factor for $(G,H)$ with respect to the generic character on
$N(F)$ determined by $\psi$ and let $\Delta_{\psi}^I$ be its extension to $G^I$ defined in Section \ref{sec:pure_tf}. We will identify the element $s \in
Z(\hat H)^\Gamma$ with its image in $\hat T_0$ under $\hat\eta$. Then from Section \ref{sec:schar_def} we have the functions $\Theta^s_\phi$ on $G^I(F)$ and
$\mc{S}\Theta_{\phi^H}$ on $H(F)$. The main result of this paper is

\begin{thm} For any strongly regular semi-simple element $\gamma \in G^I(F)$ the following equality holds
\begin{equation*}
\Theta^s_\phi(\gamma) = \sum_{\gamma^H \in H_\tx{sr}(F)/\tx{st}} \Delta_\psi^I(\gamma^H,\gamma) \frac{D(\gamma^H)^2}{D(\gamma)^2}\mc{S}\Theta_{\phi^H}(\gamma^H)
\end{equation*}
\end{thm}

In terms of the normalized $L$-packets, this statement can be reformulated as follows. Let $(\phi,u) : G \rightarrow G^u$ be a pure inner twist of the type
considered in \cite[\S4.6]{DR09} and let $\Delta_{\psi,u}$ be the normalization of the absolute transfer factor for $(G^u,H)$ corresponding to $\Delta_\psi$ as
in Section \ref{sec:pure_tf}. Then

\begin{thm} \label{thm:chid}
For any strongly regular semi-simple element $\gamma \in G^u(F)$ the following equality holds
\begin{equation} \Theta^s_{\phi,u}(\gamma) = \sum_{\gamma^H \in H_\tx{sr}(F)/\tx{st}} \Delta_{\psi,u}(\gamma^H,\gamma)
\frac{D(\gamma^H)^2}{D(\gamma)^2}\mc{S}\Theta_{\phi^H}(\gamma^H) \label{eq:chid} \end{equation}
\end{thm}


\section{Endoscopic signs} \label{sec:signs}
In this section we only need the notation from the beginning of Section \ref{sec:main}. Moreover, it is independent of the restrictions posed on $p$ in Section
\ref{sec:mainst}. The only restriction we impose on $p$ is $p>2$, although this again is just for convenience and could be removed.

There are three signs which can be assigned to the pair of groups $(G,H)$ (and some auxiliary choices) and which we need to equate. The first one is
\[ \epsilon(G,H) = (-1)^{r_G-r_H} \]
where $r_G$ and $r_H$ are the $F$-split ranks of $G$ and $H$. This sign plays an important role in the character formulas of \cite{DR09}.

The second sign enters in the normalization of the geometric transfer factors. It is defined relative to an additive character $\psi : F \rightarrow \C^\times$
as the local $\epsilon$-factor $\epsilon_L(V,\psi)$ where $V$ is the virtual representation of $\Gamma$ of degree $0$ given by the difference of the
$\Gamma$-representations $V_G:= X^*(T_0)\otimes \C$ and $V_H:= X^*(T_0^H) \otimes \C$.

The third appears in Waldspurger's work \cite{Wal95} on the local trace formula for Lie algebras. To construct it, let $\psi : F \rightarrow \C^\times$ be an
additive character and $B : \mf{g}(F) \times \mf{g}(F) \rightarrow F$ a non-degenerate, $\tx{Ad}(G(F))$-invariant, symmetric bilinear form. With this data,
Waldspurger defines in \cite[VIII]{Wal95} for a lattice $r \subset \mf{g}(F)$
\begin{eqnarray*}
I(r)&=&\int_r \psi(B(x,x)/2) dx\\
\tilde r&=&\{ x \in \mf{g}(F)| \forall y \in r\ \psi(B(x,y))=1 \}
\end{eqnarray*}
and remarks that the function
\[ r \mapsto \frac{I(r)}{|I(r)|} \]
is constant when restricted to the set $\{ r| \tilde r \subset 2r \}$. This constant he then calls $\gamma_\psi(B)$, or $\gamma_\psi(\mf{g})$ when $B$ is
understood. Furthermore, in loc. cit. Waldspurger explains how to transfer $B$ to a non-degenerate, $\tx{Ad}(H(F))$-invariant, symmetric bilinear form
$B_\mf{h}$ on $\mf{h}(F)$, thereby obtaining $\gamma_\psi(B_\mf{h})$. The second sign we are interested in is $\gamma_\psi(B)\gamma_\psi(B_\mf{h})^{-1}$. (The
word "sign" is not yet justified here, all we know is that both constants and hence their quotient are eight roots of unity. We will see however that in our
case the quotient is a sign.)

We extend the bilinear form $B$ to a symmetric bilinear form $\mf{g}(\ol{F}) \times \mf{g}(\ol{F}) \rightarrow \ol{F}$ in the obvious way and denote it by the
same letter. As remarked in loc.cit., this extension is $\tx{Ad}(G(\ol{F})) \rtimes \Gamma$-invariant. It is clear that if $V \subset \mf{g}$ is a subspace of
$\mf{g}$ defined over some extension $E$ of $F$, then the restriction of $B$ to $V$ defines a symmetric bilinear form $V(E) \times V(E) \rightarrow E$.

The purpose of this section is to prove the following
\begin{pro} \label{pro:signs}Let $\psi : F \rightarrow \C^\times$ be an additive character which is non-trivial on $O_F$ and trivial on $\pi O_F$. Let $B$ be a "good" bilinear
form in the sense of \cite[A.1]{DR09}. Then
\[ \epsilon_L(V,\psi) = \epsilon(G,H) = \gamma_\psi(B)\gamma_\psi(B_\mf{h})^{-1} \]
\end{pro}

The proof is contained in the following lemmas.

\rmk We would like to point out that the second of these equalities is also proved in \cite{KV2}. The proof given here is different from the one in loc. cit.
and establishes a connection between the above signs and the number of symmetric orbits of $\Gamma$ in $R(T^H,G)$. This number is an important invariant in
endoscopy and thus the following lemmas may be of independent interest.

\begin{lem}
\[ \epsilon(G,H) = \tx{det}(\omega) \]
\end{lem}
\pf A similar argument is given in the proof of \cite[Lemma 12.3.5]{DR09}, but we will present it here since our situation and notation are different.
$\vartheta$ is a finite-order automorphism of the real vector space $X^*(T_0) \otimes \R$ and hence is diagonalizable over $\C$ with eigenvalues roots of
unity, and all non-real eigenvalues come in conjugate pairs. Thus $\det(\vartheta)=(-1)^{\dim(V_G)-\dim(V_G^\vartheta)}$. In the same way
$\det(\omega\vartheta)=(-1)^{\dim(V_H)-\dim(V_H^{\omega\vartheta})}$. But
\[ \epsilon(G,H) = (-1)^{\dim(V_G^\Gamma)-\dim(V_H^\Gamma)}=(-1)^{\dim(V_G^\vartheta)-\dim(V_H^{\omega\vartheta})}=\det(\omega) \]
\qed

\begin{lem}
\[ \epsilon_L(V,\psi)=\det(\omega) \]
\end{lem}
\pf The $\Gamma$ representations $V_G$ and $V_H$ are unramified. Applying \cite[3.4.6]{Tat77} and noting that the isomorphism of local class field theory used
in loc. cit. is normalized so that $\tx{Fi}$ corresponds to $\pi$, we obtain
\begin{eqnarray*}
\epsilon_L(V_G-V_H,\psi)&=&\det V_G(\tx{Fi}^{-1})\det V_H(\tx{Fi}^{-1})^{-1}\\
&=&\left[ \frac{\det(\vartheta)}{\det(\omega\vartheta)} \right]^{-1}\\
&=&\det(\omega)
\end{eqnarray*}
\qed

These two lemmas complete the proof of the first equality in Proposition \ref{pro:signs}. To continue with the second equality, we need to recall some notions
from \cite{LS87}. Let $T$ be a maximal torus of $G$, and $\mc{O}$ be a $\Gamma$-orbit in $R(T,G)$, the set of roots of $T$ in $G$. Then $-\mc{O}$ is also a
$\Gamma$-orbit in $R(T,G)$ and either $\mc{O}=-\mc{O}$, in which case $\mc{O}$ is called a \ti{symmetric} orbit, or $\mc{O} \cap -\mc{O} = \emptyset$, in which
case $\mc{O}$ is called an \ti{asymmetric} orbit. For $\alpha \in R(T,G)$ let $\Gamma_\alpha$ be the stabilizer of $\alpha$ and $\Gamma_{\pm\alpha}$ be the
stabilizer of the set $\{\alpha,-\alpha\}$. Let $F_\alpha$ and $F_{\pm\alpha}$ be the fixed fields of $\Gamma_\alpha$ and $\Gamma_{\pm\alpha}$ in $\ol{F}$.
Then $[\Gamma_\alpha,\Gamma_{\pm\alpha}]$ equals $2$ if the orbit of $\alpha$ is symmetric and $1$ if it is asymmetric. If $T$ is unramified, then both
$F_\alpha$ and $F_{\pm\alpha}$ lie in $F^u$.

For any $\Gamma$-invariant subset $S \subset R(T,G)$ we put
\[ \mf{g}_S = \bigoplus_{\alpha \in S} \mf{g}_\alpha \]
This is clearly a vector subspace of $\mf{g}$ defined over $F$.

\begin{lem} Let $T$ be a maximal torus of $G$ stably conjugate to $T_0^\omega$. Then
\[ \gamma_\psi(B)\gamma_\psi(B_\mf{h})^{-1} = \prod_\mc{O} \gamma_\psi(B|_{\mf{g}_{\mc{O}}(F)}) \]
where $\mc{O}$ runs over the set of symmetric orbits of $\Gamma$ in $R(T,G)$.
\end{lem}
\pf We consider the root decomposition of $\mf{g}$ relative to $T$:
\[ \mf{g} = \mf{t} \oplus \bigoplus_{\alpha\in R(T,G)} \mf{g}_\alpha \]
If we put $\mf{g}_0=\mf{t}$ then the invariance of $B$ implies that for all $\alpha,\beta \in R(T,G) \cup \{0\}$ such that $\alpha\neq -\beta$ the subspaces
$\mf{g}_\alpha$ and $\mf{g}_\beta$ of $\mf{g}$ are orthogonal with respect to $B$. This means that if $\mc{O}_1,...,\mc{O}_k$ are the orbits in $R(T,G)$ of the
group $\Gamma \times \{\pm 1\}$, where $\{\pm 1\}$ acts by scalar multiplication, then
\[ \mf{g}(F) = \mf{t}(F) \oplus \bigoplus_{i=1}^k \mf{g}_{\mc{O}_i}(F) \]
is an orthogonal decomposition of $\mf{g}(F)$. Thus $\gamma_\psi(B)$ factors as
\[ \gamma_\psi(B) = \gamma_\psi(B|_{\mf{t}(F)}) \prod_{i=1}^k \gamma_\psi(B|_{\mf{g}_{\mc{O}_i}(F)}) \]
Consider one of the orbits $\mc{O}_i$. Either $\Gamma$ acts transitively on it, in which case it is a symmetric $\Gamma$-orbit, or it decomposes as a disjoint
union of two asymmetric $\Gamma$-orbits. We assume that the latter is the case, and write $\mc{O}_i=\mc{O}_i' \sqcup -\mc{O}_i'$ where $\mc{O}_i'$ is one of
the two $\Gamma$ orbits in $\mc{O}_i$. Then $\mf{g}_{\mc{O}_i}=\mf{g}_{\mc{O}_i'} \oplus \mf{g}_{-\mc{O}_i'}$ is a decomposition over $F$ as a direct sum of
isotropic spaces. Let $r_+ \subset \mf{g}_{\mc{O}_i'}(F)$ and $r_- \subset \mf{g}_{-\mc{O}_i'}(F)$ be large enough lattices. Then
$\gamma_\psi(B|_{\mf{g}_{\mc{O}_i}(F)})$ is by definition the complex sign of
\begin{eqnarray*}
&&\int_{r_+ \oplus r_-} \psi(B(x+y,x+y)/2) d(x,y) \\
&=&\int_{r_+} \int_{r_-} \psi(B(x,y)) dxdy
\end{eqnarray*}
For each $x \in r_+$ the map $y \mapsto \psi(B(x,y))$ is a character of the additive group $r_-$. Thus if $r_+^0$ is the subgroup of $r_+$ consisting of all
$x$ s.t. this character is trivial, the above integral is equal to the positive real constant $\vol(r_+^0,dx)\vol(r_-,dy)$. This shows
$\gamma_\psi(B|_{\mf{g}_{\mc{O}_i}(F)})=1$ and we conclude that
\[ \gamma_\psi(B) = \gamma_\psi(B|_{\mf{t}(F)}) \prod_{\mc{O}} \gamma_\psi(B|_{\mf{g}_{\mc{O}_i}(F)}) \]
where $\mc{O}$ runs over the set of symmetric $\Gamma$-orbits in $R(T,G)$.

We can apply the same reasoning to the Lie algebra $\mf{h}$ with the bilinear form $B_\mf{h}$ and the torus $T_0^H$. Since $T_0^H$ is contained in a Borel
defined over $F$, there are no symmetric orbits of $\Gamma$ in $R(T_0^H,H)$ and we conclude
\[ \gamma_\psi(B_\mf{h}) = \gamma_\psi(B_\mf{h}|_{\mf{t}_0^H(F)}) \]
But we have chosen the torus $T$ so that there exists an admissible isomorphism $T_0^H \rightarrow T$ over $F$, and the bilinear form $B_\mf{h}$ is constructed
so that the differential of this admissible isomorphism identifies $B_\mf{h}|_{\mf{t}_0^H(F)}$ with $B|_{\mf{t}(F)}$. Thus
\[ \gamma_\psi(B_\mf{h}) = \gamma_\psi(B|_{\mf{t}(F)}) \]
and the lemma now follows. \qed

\begin{lem} Let $\mc{O}$ be a symmetric orbit of $\Gamma$ in $R(T,G)$. Then
\[ \gamma_\psi(B|_{\mf{g}_\mc{O}(F)}) = -1 \]
\end{lem}
\pf Choose $\alpha \in \mc{O}$ and $\sigma_\alpha \in \Gamma_{\pm\alpha} \smallsetminus \Gamma_\alpha$. We can choose a non-zero $E \in \mf{g}_\alpha(F_\alpha)
\cap [\mf{g}_o \smallsetminus \mf{g}_{o^+}]$ and then we have $\sigma_\alpha(E) \in \mf{g}_{-\alpha}(F_\alpha) \cap [\mf{g}_o \smallsetminus \mf{g}_{o^+}]$.
Then by \cite[\S A.1]{DR09}
\[B(E,\sigma_\alpha(E)) \in O_{F_{\pm\alpha}}^\times \]
The map
\[ \phi : F_\alpha \rightarrow \mf{g}_\mc{O}(F),\qquad \lambda \mapsto \sum_{\sigma \in \Gamma/\Gamma_\alpha} \sigma(\lambda E) \]
is an isomorphism of $F$-vector spaces and clearly $\gamma_\psi(B|_{\mf{g}_\mc{O}(F)})=\gamma_\psi(\phi^*B)$. To compute the bilinear form $\phi^*B : F_\alpha
\times F_\alpha \rightarrow F$ we notice that if $\sigma_1,...,\sigma_k$ are representatives for $\Gamma/\Gamma_{\pm\alpha}$, then
\[ \mf{g}_\mc{O} = \bigoplus_{i=1}^k (\mf{g}_{\sigma_i(\alpha)} \oplus \mf{g}_{\sigma_i(-\alpha)}) \]
is an orthogonal sum of hyperbolic planes. Then a direct computation shows that
\[ \phi^*B(\lambda,\mu) = \tr_{F_{\pm\alpha}/F}\Big( [\lambda\sigma_\alpha(\mu)+\mu\sigma_\alpha(\lambda)]B(E,\sigma_\alpha(E)) \Big) \]
If we put
\begin{eqnarray*}
\psi'(x)&=&\tr_{F_{\pm\alpha}/F}(B(E,\sigma_\alpha(E))x)\\
B'(\mu,\lambda)&=&\lambda\sigma_\alpha(\mu)+\mu\sigma_\alpha(\lambda)
\end{eqnarray*}
then $\psi' : F_{\pm\alpha} \rightarrow \C^\times$ is an additive character and $B' : F_\alpha \times F_\alpha \rightarrow F_{\pm\alpha}$ is a non-degenerate
$F_{\pm\alpha}$-bilinear form, and clearly $\gamma_\psi(\phi^*B) = \gamma_{\psi'}(B')$.

We will now compute $\gamma_{\psi'}(B')$.

First we claim that $\psi'$ is non-trivial on $O_{F_{\pm\alpha}}$ but trivial on $\pi O_{F_{\pm\alpha}}$. To see this, note that $\tr_{F_{\pm\alpha}/F}$
induces for each $i \in \Z$ a homomorphism of additive groups $\pi^i O_{F_{\pm\alpha}} \rightarrow \pi^i O_F$ which fits into the diagram
\begin{diagram}
\pi^i O_{F_{\pm\alpha}}&\rTo\pi^i& O_F\\
\dTo<{\small\textrm{mod } \pi^{i+1}}&&\dTo>{\small \textrm{mod } \pi^{i+1}}\\
k_{F_{\pm\alpha}}&\rTo^{\tr}&k_F
\end{diagram}
and thus $\tr_{F_{\pm\alpha}/F} : O_{F_{\pm\alpha}} \rightarrow O_F$ is surjective (\cite[V.\S1.Lemma 2]{Ser79}). This together with $B(E,\sigma_\alpha(E)) \in
O_{F_{\pm\alpha}}^\times$ implies the claim about $\psi'$.

Next we compute the dual of the $O_{F_{\pm\alpha}}$-lattice $O_{F_\alpha}$ with respect to $\psi'\circ B'$.
\begin{eqnarray*}
&&\{ x \in F_\alpha|\ \forall y \in O_{F_{\alpha}}: \psi'(B(x,y))=1 \}\\
&=&\{ x \in F_\alpha|\ \forall y \in O_{F_{\alpha}}: B(x,y) \in \pi O_{F_{\pm\alpha}} \}\\
&=&\pi\{ x \in F_\alpha|\ \forall y \in O_{F_{\alpha}}: x\sigma_\alpha(y)+y\sigma_\alpha(x) \in O_{F_{\pm\alpha}} \}\\
&=&\pi\{ x \in F_\alpha|\ \forall y \in O_{F_{\alpha}}: xy+\sigma_\alpha(y)\sigma_\alpha(x) \in O_{F_{\pm\alpha}} \}
\end{eqnarray*}
Thus we are looking for $\pi$ times the dual of $O_{F_\alpha}$ with respect to the bilinear form $(x,y) \mapsto \tr_{F_\alpha/F_{\pm\alpha}}(xy)$. This dual is
the codifferent of $F_\alpha/F_{\pm\alpha}$, which equals $O_{F_\alpha}$ since $F_\alpha/F_{\pm\alpha}$ is an unramified extension.

We conclude that the lattice $O_{F_\alpha}$ has the property that it contains its dual with respect to $\psi'\circ B'$. Since we are imposing the restriction
$p>2$ and thus $O_{F_\alpha}=2O_{F_\alpha}$. Then by definition, $\gamma_{\psi'}(B')$ is the complex sign of
\[ I := \int_{O_{F_\alpha}} \psi'(N(x)) dx \]
where $N : F_\alpha \rightarrow F_{\pm\alpha}$ is the norm map and $dx$ is a Haar measure on the additive group $F_\alpha$. Let $(\xi_k)_{k \in k_{F_\alpha}}$
be a set of representatives for $O_{F_\alpha}/\pi O_{F_\alpha}$. Then
\[ I = \sum_{k \in k_{F_\alpha}} \int_{\pi O_{F_\alpha}} \psi'(N(\xi_k+x)) dx \]
One computes immediately that $\psi'(N(\xi_k+x))=\psi'(N(\xi_k))$ for all $x \in \pi O_{F_\alpha}$ since $\psi'$ is trivial on $\pi O_{F_{\pm\alpha}}$. This
leads to
\[ I = \vol(\pi O_{F_\alpha},dx)\sum_{k \in k_{F_\alpha}} \psi'(N(\xi_k)) \]
The restriction of $\psi'$ to $O_{F_{\pm\alpha}}$ factors through the natural projection $O_{F_{\pm\alpha}} \rightarrow k_{F_{\pm\alpha}}$, and the composition
of $N$ with this projection factors through the projection $O_{F_\alpha} \rightarrow k_{F_\alpha}$ and induces the norm map associated to the extension
$k_{F_\alpha}/k_{F_{\pm\alpha}}$, which we also call $N$. Thus
\begin{eqnarray*}
I&=&\vol(\pi O_{F_\alpha},dx)\sum_{k \in k_{F_\alpha}} \psi'(N(k)) \\
&=&\vol(\pi O_{F_\alpha},dx)\left[ 1+ \sum_{k \in k_{F_\alpha}^\times} \psi'(N(k)) \right]
\end{eqnarray*}
Now $N: k_{F_\alpha}^\times \rightarrow k_{F_{\pm\alpha}}^\times$ is a surjective homomorphism, the cardinality of whose fibers we will call $A$. Then
\begin{eqnarray*}
I&=&\vol(\pi O_{F_\alpha},dx)\left[ 1+ A\sum_{k \in k_{F_{\pm\alpha}}^\times} \psi'(k) \right]\\
&=&\vol(\pi O_{F_\alpha},dx)\left[ -(A-1)+ A\sum_{k \in k_{F_{\pm\alpha}}} \psi'(k) \right]\\
&=&-(A-1)\vol(\pi O_{F_\alpha},dx)
\end{eqnarray*}
since $\psi'$ is a non-trivial character on the additive group $k_{F_{\pm\alpha}}$. We conclude that $I$ is a negative real number, and the lemma follows.\qed

\begin{lem} \[ \det(\omega)=(-1)^N \]
where $N$ is the number of symmetric orbits of $\Gamma$ in $R(T,G)$.
\end{lem}
\pf We choose a $g \in G(\ol{F})$ s.t. $\tx{Ad}(g) : T_0^\omega \rightarrow T$ is an isomorphism over $F$ and use it to regard $\omega$ and $\vartheta$ as
automorphisms of $R(T,G)$. Moreover put $B=\tx{Ad}(g)B_0$ and write $\alpha>0$ if $\alpha \in R(T,B)$. Let
\begin{eqnarray*}
S&=&\{ \alpha \in R(T,G)| \alpha >0 \wedge \omega\alpha < 0 \} \\
S'&=&\{ \alpha \in R(T,G)| \alpha >0 \wedge \omega\vartheta\alpha < 0 \}
\end{eqnarray*}
Since $\vartheta$ preserves the set of positive roots in $R(T,G)$, it induces a bijection $S' \rightarrow S$. Thus
\[ \det(\omega)=(-1)^{|S|} = (-1)^{|S'|} \]

\ul{Claim 1:} The cardinality of $S'$ is congruent mod $2$ to the cardinality of the intersection of $S'$ with the union of the symmetric orbits of $\Gamma$ in
$R(T,G)$.

Put $T=\omega\vartheta$ for short. Then $\Gamma$ acts on $R(T,G)$ via the cyclic group $<T>$. Let $\mc{O}$ be an orbit. We claim that the sets
\begin{eqnarray*}
\mc{O}_+&=&\{\alpha \in \mc{O}|\ \alpha>0 \wedge T\alpha < 0 \}\\
\mc{O}_-&=&\{\alpha \in \mc{O}|\ \alpha<0 \wedge T\alpha > 0 \}\\
\end{eqnarray*}
have the same cardinality. To see this, consider the directed graph in the vector space $X^*(T_\tx{ad})\otimes \R$ whose vertices are given by $\mc{O}$ and
whose edges are given by
\[ \{ (\alpha,T\alpha)|\ \alpha \in \mc{O} \} \]
Then $\mc{O}_+$ is in bijection with the set of edges which start in the positive half space of $X^*(T_\tx{ad})\otimes \R$ and end in the negative, while
$\mc{O}_-$ is in bijection with the set of edges which start in the negative half space and end in the positive. But our graph is a closed loop, so these sets
must have the same cardinality.

If $\mc{O}$ is an asymmetric orbit, then $-\mc{O}$ is also one and is disjoint from $\mc{O}$, and multiplication by $-1$ gives a bijection $\mc{O}_-
\rightarrow [-\mc{O}]_+$. We conclude that $S' \cap (\mc{O} \cup -\mc{O})$ has an even cardinality. This proves Claim 1.

\ul{Claim 2:} Let $\mc{O}$ be a symmetric orbit. Then its intersection with $S'$ has an odd cardinality.

The group $<T>$ acts on $\mc{O}/\{\pm 1\}$ and all elements of the latter set are of the form $\{\alpha,-\alpha\}$ with $\alpha \in \mc{O}$. We choose an
element $A \in \mc{O}/\{\pm 1\}$, and let $n=|\mc{O}|/2-1$. Then $A,TA,...,T^nA$ enumerates $\mc{O}/\{\pm 1\}$. For each $0 \leq i \leq n$ let $\alpha_i$ be
the positive member of $T^iA$. Then for each such $i$ one of two cases occurs: either $T\alpha_i=-\alpha_{i+1}$ and $T(-\alpha_i)=\alpha_{i+1}$, or
$T\alpha_i=\alpha_{i+1}$ and $T(-\alpha_i)=-\alpha_{i+1}$ (where we adopt the convention $\alpha_{n+1}=\alpha_0$). The cardinality of $S' \cap \mc{O}$ is the
number of $0 \leq i \leq n$ for which the first case occurs. Now let $M$ be the number of $0 \leq i < n$ for which the first case occurs (note the sharp
inequality!). If $M$ is even, then $T^n\alpha_0=\alpha_n$ and thus $T\alpha_n$ must equal $-\alpha_0$, for otherwise the set
$\{\alpha_0,T\alpha_0,...,T^n\alpha_0\}$ will be a $T$-invariant subset of $\mc{O}$, which is impossible. Thus $|S'\cap\mc{O}|=M+1$ is an odd number. If
conversely $M$ is odd, then $T^n\alpha_0=-\alpha_n$ and by the same reasoning $T(-\alpha_n)=(-\alpha_0)$. It follows then that $|S'\cap\mc{O}|=M$, again an odd
number. This proves Claim 2.

The two claims together imply that $(-1)^{|S'|}=(-1)^N$ and this finishes the lemma.\qed

The second equality in Proposition \ref{pro:signs} now follows from these lemmas.


\section{A formula for the unstable character} \label{sec:schar}

The purpose of this section is to establish a reduction formula, similar to the ones in \cite[\S9,\S10]{DR09}, for $\Theta^{t}_{\rho,u}$. Before we can do so,
we need some cohomological facts.


\subsection{Cohomological lemmas I}
We begin by recalling some well-known basic facts about Tate-Nakayama duality as used in endoscopy. For this, we will deviate from the notation established so
far in order to make the statements in their natural generality. Let $E/F$ be a finite extension of local fields of characteristic 0, $\Gamma = \tx{Gal}(E/F)$,
$u_{E/F} \in H^2(\Gamma,E^\times)$ the canonical class of the extension $E/F$, $T$ a torus over $F$ which splits over $E$, and $\hat T$ its dual complex torus.

\begin{lem} \label{lem:tndseq}
We have the exact sequences
\begin{diagram}
1&\rTo&(\hat T^\Gamma)^\circ&\rTo&\hat T^\Gamma&\rTo&H^1(\Gamma,X_*(\hat T))&\rTo&0\\
0&\rTo&X^*(\hat T/\hat T^\Gamma)&\rTo&X^*(\hat T/(\hat T^\Gamma)^\circ)&\rTo&H_T^{-1}(\Gamma,X^*(\hat T))&\rTo&0
\end{diagram}
\end{lem}
\pf For the first one, tensor the exponential sequence
\[ 0 \rightarrow \Z \rightarrow \C \stackrel{e^{2\pi iz}}{\longrightarrow} \C^\times \rightarrow 1 \]
with $X_*(T)$ and take $\Gamma$-invariants, noting that the image of $[X_*(\hat T)\otimes \C]^\Gamma=\tx{Lie}(\hat T)^\Gamma$ in $\hat T^\Gamma$ under the
exponential map is $(\hat T^\Gamma)^\circ$.

For the second one, observe that an element of $X^*(\hat T)$ is in the kernel of the norm map precisely when it is trivial on $(\hat T^\Gamma)^\circ$ and in
the augmentation submodule precisely when it is trivial on $\hat T^\Gamma$.\qed

\begin{lem} \label{lem:tndpairings} The following three pairings
\[ H^1(\Gamma,X_*(\hat T)) \otimes H_T^{-1}(\Gamma,X^*(\hat T)) \rightarrow \C^\times \]
are equal.
\begin{enumerate}
\item The pairing induced by the standard pairing $\hat T \times X^*(\hat T) \rightarrow \C^\times$ via the above sequences.
\item
\begin{diagram}[height=16pt,width=30pt,tight,midshaft,scriptlabels]
H^1(\Gamma,X_*(\hat T))\\
\otimes&&\rTo^{\cup}&H_T^0(\Gamma,\Z)&\rEquals&\Z/|\Gamma|\Z&\rTo^{\cdot|\Gamma|^{-1}}&\Q/\Z&\rTo^{e^{2\pi iz}}&\C^\times\\
H_T^{-1}(\Gamma,X^*(\hat  T))
\end{diagram}
\item
\begin{diagram}[height=16pt,width=30pt,tight,midshaft,scriptlabels]
H^1(\Gamma,X_*(\hat T))&&\rEquals&H^1(\Gamma,X^*(T))\\
&&&\otimes&\rTo^{\cup}&H^2(\Gamma,E^\times)&\rTo^{\tx{inv}}&\Q/\Z&\rTo^{e^{2\pi i z}}&\C^\times\\
H_T^{-1}(\Gamma,X_*(T))&&\rTo^{\cup u_{E/F}}&H^1(\Gamma,T)\\
\end{diagram}
\end{enumerate}
\end{lem}

\pf The equality of the pairings in 2. and 3. is an immediate consequence of local class field theory, more precisely of the following commutative square.
\begin{diagram}[midshaft]
H_T^2(\Gamma,E^\times)&\lTo^{\cup u_{E/F}}&&H_T^0(\Gamma,\Z)\\
\dTo<{\tx{inv}}&&&\dEquals\\
|\Gamma|^{-1}\Z/\Z&\lTo^{\cdot|\Gamma|^{-1}}&&\Z/|\Gamma|\Z
\end{diagram}
In order to relate pairings 1. and 2. take $t \in \hat T^\Gamma$ and $\phi \in X^*(\hat T/(\hat T^\Gamma)^\circ)$. Choose $z \in \tx{Lie}(\hat T)=X_*(\hat
T)\otimes \C$ mapping to $t$ under the exponential map. Then the image of $t$ in $H^1(\Gamma,X_*(\hat T))$ is represented by the cocycle $\tau \mapsto \tau
z-z$. Now using the appropriate cup product formula and denoting the canonical pairing $X^*(\hat T) \otimes X_*(\hat T) \rightarrow \Z$ by $\langle\rangle$ we
compute
\begin{eqnarray*}
(\tau z-z) \cup \phi&=&\sum_{\tau \in \Gamma}\langle \tau\phi,\tau z-z\rangle\\
&=&|\Gamma|\langle \phi,z \rangle
\end{eqnarray*}
Note that we have used that $\phi$ is in the kernel of the norm map. It follows that
\[ \exp(2\pi i |\Gamma|^{-1}(\tau z-z)\cup \phi) = \exp(2\pi i \langle \phi,z \rangle) = \langle \phi,t \rangle \]
\qed

\begin{lem} \label{lem:tndiso}
Assume that $E/F$ is either an unramified extension of $p$-adic fields, or $\C/\R$. In the $p$-adic case, let $\pi \in E^\times$ be a uniformizer and $\sigma
\in \Gamma$ be the Frobenius element. In the real case, let $\pi=-1$ and $\sigma \in \Gamma$ be complex conjugation. Then the map
\[ [\lambda] \mapsto \lambda(\pi) \]
induces the same isomorphism
\[ [X_*(T)_\Gamma]_\tx{tor} = H^{-1}_T(\Gamma,X_*(T)) \rightarrow H^1(\Gamma,T) \]
as the isomorphism given by $\cup u_{E/F}$. Here we regard $\lambda(\pi) \in T(E)$ as the class in $H^1(\Gamma,T)$ represented by the unique element $z \in
Z^1(\Gamma,T)$ s.t. $z(\sigma)=\lambda(\pi)$.
\end{lem}
\pf By definition, $H^{-1}_T(\Gamma,X_*(T)) = \tx{Ker}( N : X_*(T) \rightarrow X_*(T) )/IX_*(T)$ where $N$ is the norm map and $I \subset \Z[\Gamma]$ is the
augmentation ideal. If $\lambda \in X_*(T)$ is torsion modulo $IX_*(T)$, then some multiple of it is killed by $N$, and since $X_*(T)$ is torsion-free this
means that $\lambda$ itself is killed by $N$. Thus
\[[X_*(T)_\Gamma]_\tx{tor} \subset H^{-1}_T(\Gamma,X_*(T)) \]
The converse inclusion follows from the finiteness of $H^{-1}_T(\Gamma,X_*(T))$. This justifies the first equality.

It is well known from local class field theory that the fundamental class of $E/F$ is represented by the 2-cocycle
\[ (\sigma^a,\sigma^b) \mapsto \begin{cases} 1&, 0\leq a+b < |\Gamma|\\ \pi&,\tx{ else } \end{cases} \]
If $\lambda \in X_*(T)$ is torsion modulo $IX_*(T)$, then applying the appropriate cup-product formula one sees
\[ ([\lambda] \cup u_{E/F})(\sigma) = \sum_{i=0}^{|\Gamma|-1} \sigma^{i+1}\lambda( \pi^{\tx{char}_{\{i+1\geq |\Gamma|\}}} ) = \lambda(\pi) \]
\qed

This isomorphism is sometimes called the Tate-Nakayama isomorphism. We will denote it by $\tx{TN}$. In the case that $E/F$ is an unramified extension of
$p$-adic fields, DeBacker and Reeder construct in \cite[Cor 2.4.3]{DR09} another isomorphism
\[ [X_*(T)_\Gamma]_\tx{tor} \rightarrow H^1(\Gamma,T) \]
We will call this isomorphism $\tx{DR}$. It turns out that these two isomorphisms are almost identical, namely

\begin{lem} \label{lem:tndanticom} The following diagram commutes
\begin{diagram}
[X_*(T)_\Gamma]_\tx{tor}&&\rTo^{\lambda\mapsto-\lambda}&&[X_*(T)_\Gamma]_\tx{tor} \\
&\rdTo<{\tx{DR}}&&\ldTo>{\tx{TN}}&\\
&&H^1(\Gamma,T)&&
\end{diagram}
\end{lem}
\pf By construction, $\tx{DR}$ sends $[\lambda] \in [X_*(T)_\Gamma]_\tx{tor}$ to the class in $H^1(\Gamma,T)$ of the unique cocycle $z$ whose value at
$\tx{Fi}$ equals $t_\lambda=\lambda(\pi)$, while $\tx{TN}$ sends $[\lambda]$ to the class in $H^1(\Gamma,T)$ of the unique cocycle $z'$ whose value at
$\sigma=\tx{Fi}^{-1}$ equals $\lambda(\pi)$. But
\[z(\sigma)=\sigma(z(\tx{Fi})^{-1})=\sigma(\lambda(\pi)^{-1})=\sigma(\lambda)^{-1}(\pi) \]
Since $\lambda$ and $\sigma(\lambda)$ give rise to the same element of $X_*(T)_\Gamma$, the lemma follows.\qed


\subsection{Cohomological lemmas II}
We now return to the previously established notation. Recall the diagram \eqref{diag}.

We call $a_G$ the composition
\[ H^1(\Gamma,G) \rightarrow \tx{Irr}(\pi_0(Z(\hat G)^\Gamma)) \]
of the right vertical isomorphisms in this diagram. In \cite[Thm 1.2]{Kot86} Kottwitz defines another isomorphism
\[ H^1(\Gamma,G) \rightarrow \tx{Irr}(\pi_0(Z(\hat G)^\Gamma)) \]
which he calls $\alpha_G$. This isomorphism can be normalized in two different ways, and the two normalizations differ by a sign.

\begin{lem} \label{lem:kotdr}
Depending on the normalization of $\alpha_G$, one has
\[ a_G = \pm\alpha_G \]
\end{lem}
\pf Assume first that $G=T$ is a torus. One normalization of the isomorphism $\alpha_G$ is then given by the composition
\[ H^1(\Gamma,T) \rightarrow \tx{Irr}(H^1(\Gamma,X^*(T))) \rightarrow \tx{Irr}(\pi_0(\hat T^\Gamma)) \]
where the first map arises via the cup product pairing
\[  H^1(\Gamma,X^*(T)) \otimes H^1(\Gamma,T) \rightarrow \C^\times \]
and the second map is the dual of the isomorphism $\pi_0(\hat T^\Gamma) \rightarrow H^1(\Gamma,X_*(\hat T))$ of Lemma \ref{lem:tndseq}.

Thus if we precompose $\alpha_G$ by $\tx{TN}$ then by Lemma \ref{lem:tndpairings} the resulting isomorphism
\[ [X_*(T)_\Gamma]_\tx{tor} \rightarrow \tx{Irr}(\pi_0(\hat T^\Gamma)) \]
will be given by the standard pairing $\hat T \times X^*(\hat T) \rightarrow \C^\times$.

On the other hand if we precompose $a_G$ by $\tx{TN}$ then by Lemma \ref{lem:tndanticom} the resulting isomorphism
\[ [X_*(T)_\Gamma]_\tx{tor} \rightarrow \tx{Irr}(\pi_0(\hat T^\Gamma)) \]
will be given by the negative of the standard pairing $\hat T \times X^*(\hat T) \rightarrow \C^\times$.

This proves that in the case $G=T$ with our normalization of $\alpha_G$ we have $a_G=-\alpha_G$. For the general case let $T \subset G$ be an elliptic maximal
torus and consider the commutative diagrams

\begin{diagram}
\tx{Irr}(\pi_0(\hat T^\Gamma))&\rTo&\tx{Irr}(\pi_0(Z(\hat G)^\Gamma))  &&&  \tx{Irr}(\pi_0(\hat T^\Gamma))&\rTo&\tx{Irr}(\pi_0(Z(\hat G)^\Gamma))\\
\uTo<{a_T}&&\uTo>{a_G}  &&&  \uTo<{-\alpha_T}&&\uTo>{-\alpha_G}\\
H^1(\Gamma,T)&\rTo&H^1(\Gamma,G) &&& H^1(\Gamma,T)&\rTo&H^1(\Gamma,G)
\end{diagram}

The fact that the right diagram commutes is part of the statement of \cite[Thm 1.2]{Kot86}, while for the left diagram it follows from the construction. We
just proved that the left vertical arrows in the two diagrams coincide. But since $T$ is elliptic, the bottom horizontal maps are surjective by \cite[Lemma
10.1]{Kot86}. Thus the right vertical maps in the two diagrams must also coincide. \qed

\begin{lem} \label{lem:list}
Let $Q_0 \in \tx{Lie}(S_0)(F)$ be a regular semi-simple element and $\lambda \in r^{-1}(u)$. Put $Q_\lambda := \tx{Ad}(q_\lambda q_0^{-1})Q_0$. Then
\begin{enumerate}
\item  $Q_\lambda \in \tx{Lie}(S_\lambda)(F)$.
\item The image of $\lambda$ under the map
\[ [X_\Gamma]_\tx{tor} \stackrel{\tx{DR}}{\longrightarrow} H^1(\Gamma,T_0^w) \stackrel{\tx{Ad}(q_0)}{\longrightarrow} H^1(\Gamma,S_0) \]
equals $\tx{inv}(Q_0,Q_\lambda)$.
\item The map $\lambda \mapsto Q_\lambda$ establishes a bijection from $[r^{-1}(u)]$ to a set of representatives for the conjugacy classes of elements in
    $\tx{Lie}(G^u)(F)$ stably conjugate to $Q_0$.
\item Let $t \in [\hat T_0^w]^\Gamma$ and $t_{q_0}$ be its image under the dual of $\tx{Ad}(q_0^{-1}) : S_0 \rightarrow T_0^w$. Then
\[ \langle \tx{inv}(Q_0,Q_\lambda),t_{q_0} \rangle^{-1} = \lambda(t) \]
\end{enumerate}
\end{lem}
\pf Recall from \cite[\S2.8]{DR09} the equations
\[ q_\lambda^{-1}u\tx{Fi}_G(q_\lambda) =t_\lambda \dot w \qquad \qquad q_0^{-1}\tx{Fi}_G(q_0) = \dot w \]
where $t_\lambda=\lambda(\pi)$. The inner twist $\psi$ is unramified, so $Q_\lambda \in \tx{Lie}(S_\lambda)(F^u)$. To prove that $Q_\lambda \in
\tx{Lie}(S_\lambda)(F)$ it is enough to show that it is fixed by $\tx{Fi_{G^u}}=\tx{Ad}(u)\circ\tx{Fi_G}$.
\begin{eqnarray*}
\tx{Ad}(u)\tx{Fi}_G(Q_\lambda)&=&\tx{Ad}(u)\tx{Fi}_G\tx{Ad}(q_\lambda q_0^{-1})Q_0\\
&=&\tx{Ad}(u\tx{Fi}_G(q_\lambda q_0^{-1}))Q_0\\
&=&\tx{Ad}(q_\lambda\dot w \tx{Fi}_G(q_0^{-1})Q_0\\
&=&\tx{Ad}(q_\lambda\dot w\dot w^{-1}q_0^{-1})Q_0\\
&=&Q_\lambda
\end{eqnarray*}
This proves the first assertion.

By construction the element $\tx{inv}(Q_0,Q_\lambda)$ is given by the cocycle
\[ \sigma \mapsto q_0q_\lambda^{-1}u\sigma(q_\lambda q_0^{-1}) \]
We compute the value of this cocycle at $\tx{Fi}$
\begin{eqnarray*}
q_0q_\lambda^{-1}u\tx{Fi}_G(q_\lambda q_0^{-1})&=&q_0t_\lambda\dot w \tx{Fi}_G(q_0^{-1})\\
&=&\tx{Ad}(q_0)(t_\lambda\dot w (q_0^{-1}\tx{Fi}_G(q_0))^{-1})\\
&=&\tx{Ad}(q_0)(t_\lambda)
\end{eqnarray*}
This proves the second assertion.

The third assertion follows immediately from the second and Lemma \ref{lem:stcit} (or rather from its Lie-algebra analog, which is proved in the exact same
way).

Finally, by functoriality of the Tate-Nakayama pairing we have
\[ \langle \tx{inv}(Q_0,Q_\lambda), t_{q_0} \rangle^{-1} = \langle \tx{Ad}(q_0^{-1})\tx{inv}(Q_0,Q_\lambda),t \rangle^{-1} \]
By the second assertion and Lemma \ref{lem:tndanticom} the element $\tx{Ad}(q_0^{-1})(\tx{inv}(Q_0,Q_\lambda))^{-1}$ of $H^1(\Gamma,T_0^w)$ is the image of
$\lambda$ under the Tate-Nakayama isomorphism
\[ H^{-1}_T(\Gamma,X_*(T_0^w)) \rightarrow H^1(\Gamma,T_0^w) \]
and hence by Lemma \ref{lem:tndpairings} we have
\[ \langle \tx{Ad}(q_0^{-1})\tx{inv}(Q_0,Q_\lambda)^{-1},t \rangle = \lambda(t) \]

\qed


\subsection{A reduction formula for the unstable character}
We now return to the computation of $\Theta^t_{\rho,u}$.

The map
\[ [X_\Gamma]_\tx{tor} \rightarrow \tx{Irr}(C_\phi),\qquad \lambda \mapsto \rho_\lambda \]
identifies $[r^{-1}(u)]$ with $\tx{Irr}(C_\phi,u)$. Since it is given simply by restriction of characters, we have $ \tr\rho_\lambda(t) = \lambda(t)$. Moreover
$e(G^u)=\epsilon(G,G^u)$, so
\[ \Theta^t_{\rho,u} = \epsilon(G,G^u)\sum_{\lambda \in [r^{-1}(u)]} \lambda(t) \Theta_{\pi_u(\phi,\rho_\lambda)} \]

Our first goal is to use the results of \cite[\S9,\S10]{DR09} to derive a formula for $\Theta_{\pi_u(\phi,\rho_\lambda)}$ which is suitable for our purposes.
Recall that there is a depth-zero character $\theta : T_0^w(F) \rightarrow \C^\times$ determined by the Langlands parameter $\phi$.

\begin{lem} \label{lem:onechar}
Let $\lambda \in r^{-1}(u)$, $\theta_\lambda = \tx{Ad}(q_\lambda)_*\theta$, and $Q_\lambda \in \tx{Lie}(S_\lambda)(F)$ be any fixed regular semi-simple
element. Then for any $\gamma \in G^{u}_\tx{sr}(F)_0$ and any $z \in Z(F)$ we have
\[ \Theta_{\pi_u(\phi,\rho_\lambda)}(z\gamma) = \epsilon(G^{u},A_G)\theta(z)\sum_{Q} R(G^{u}_{\gamma_s},S_Q,1)(\gamma_u)[\phi_{Q_\lambda,Q}]_*\theta_\lambda(\gamma_s) \]
where $S_Q = \tx{Cent}(Q,G^{u})$ and the sum runs over any set of representatives for the $G^{u}_{\gamma_s}(F)$-conjugacy classes inside the
$G^{u}(F)$-conjugacy class of $Q_\lambda$.
\end{lem}
\pf By \cite[Lemmas 9.3.1,9.6.2]{DR09} we know
\[ \Theta_{\pi_u(\phi,\rho_\lambda)}(z\gamma) = \epsilon(G^{u},S_\lambda)\theta(z)R(G,S_\lambda,\theta_\lambda)(\gamma) \]
We will apply \cite[Lemma 10.0.4]{DR09} to the last factor, but first we want to study the indexing set of the sum appearing in the formula of that lemma. This
indexing set is
\[ Y:= \{ (S',\theta') \in \tx{Ad}(G^{u}(F))(S_\lambda,\theta_\lambda)|\ \gamma_s \in S' \} / \tx{Ad}(G^{u}_{\gamma_s}(F)) \]
First we claim that the map
\begin{eqnarray*}
\tx{Ad}(G^{u}(F))Q_\lambda&\rightarrow&\tx{Ad}(G^{u}(F))(S_\lambda,\theta_\lambda)\\
Q&\mapsto&\phi_{Q_\lambda,Q}(S_\lambda,\theta_\lambda)
\end{eqnarray*}
is a bijection. It is clearly well-defined, and is moreover surjective because if $\tx{Ad}(g)(S_\lambda,\theta_\lambda)$ belongs to the right hand side, then
$\tx{Ad}(g)Q_\lambda$ belongs to the left hand side and is a preimage. For the injectivity let
\[(S',\theta') = \phi_{Q_\lambda,Q}(S_\lambda,\theta_\lambda) = \phi_{Q_\lambda,Q'}(S_\lambda,\theta_\lambda) \]
Then $\phi_{Q',Q} \in \Omega(S',G^{u})$ and $\phi_{Q',Q}\theta' = \theta'$. Since $\theta'$ is regular, $\phi_{Q',Q}=1$ and thus $Q'=Q$.

This proves the claimed bijectivity. Moreover, since $\phi_{Q_\lambda,Q}(S_\lambda) =S_Q$ and
\[ \gamma_s \in S_Q \Rightarrow S_Q \subset G_{\gamma_s} \Rightarrow Q \in \tx{Lie}(S_Q) \subset \tx{Lie}(G^{u}_{\gamma_s}) \Rightarrow \gamma_s \in S_Q \]
we see that our bijection restricts to the bijection
\begin{eqnarray*}
\tx{Ad}(G^{u}(F))Q_\lambda \cap \tx{Lie}(G^{u}_{\gamma_s})(F)&\rightarrow&\{ (S',\theta') \in \tx{Ad}(G^{u}(F))(S_\lambda,\theta_\lambda)|\ \gamma_s \in S'\}\\
Q&\mapsto&\phi_{Q_\lambda,Q}(S_\lambda,\theta_\lambda)
\end{eqnarray*}
Both sides of this bijection carry a natural action of $G^{u}_{\gamma_s}(F)$ and that the bijection is equivariant with respect to these actions. Thus if we
put
\[ Y':=[\tx{Ad}(G^{u}(F))Q_\lambda \cap \tx{Lie}(G^{u}_{\gamma_s})(F)]/ \tx{Ad}(G^{u}_{\gamma_s}(F)) \]
we obtain a bijection
\[ Y' \rightarrow Y \]
Applying now \cite[Lemma 10.0.4]{DR09} we obtain
\[ \Theta_{\rho_\lambda}(z\gamma) = \epsilon(G^{u},S_\lambda)\theta(z) \sum_{[Q] \in Y'}
    R(G^{u}_{\gamma_s},S_Q,1)(\gamma_u)[\phi_{Q_\lambda,Q}]_*\theta_\lambda(\gamma_s) \]
To complete the lemma, we only need to observe that since $S_\lambda/Z$ is anisotropic, the maximal split subtorus of $S_\lambda$ is $A_G$ and thus
$\epsilon(G^{u},S_\lambda)=\epsilon(G^{u},A_G)$.\qed

We are now ready to establish the reduction formula for the $t$-unstable character.

\begin{pro} \label{pro:schar}
Let $Q_0 \in \tx{Lie}(S_0)(F)$ be a regular semi-simple element, $\theta_0 = \tx{Ad}(q_0)_*\theta$, and $t_{q_0}$ be the image of $t$ under the dual of
$\tx{Ad}(q_0^{-1})$. Then for any $\gamma \in G^{u}_\tx{sr}(F)_0$ and any $z \in Z(F)$ we have
\begin{eqnarray*}
\Theta^t_{\rho,u}(z\gamma)=\epsilon(G,A_G)\theta(z)&\cdot&\sum_P [\phi_{Q_0,P}]_*\theta_0(\gamma_s) \\
&&\sum_Q \langle \tx{inv}(Q_0,Q),t_{q_0}\rangle^{-1}R(G^u_{\gamma_s},S_Q,1)(\gamma_u)
\end{eqnarray*}
where $P$ runs over a set of representatives for the $G^u_{\gamma_s}$-stable classes of elements of $\tx{Lie}(G^u_{\gamma_s})(F)$ which are $G^u$-stable
conjugate to $Q_0$, and $Q$ runs over a set of representatives for the $G^u_{\gamma_s}(F)$-conjugacy classes inside the $G^u_{\gamma_s}$-stable class of $P$.
\end{pro}
\pf For each $\lambda \in r^{-1}(u)$ put $\theta_\lambda = \tx{Ad}(q_\lambda)_*\theta$ and $Q_\lambda = \tx{Ad}(q_\lambda q_0^{-1})Q_0$. Then by Lemma
\ref{lem:list} we know that $Q_\lambda \in \tx{Lie}(S_\lambda)(F)$ is regular semi-simple, and so applying Lemma \ref{lem:onechar} and using the transitivity
of the sign $\epsilon(\cdot,\cdot)$ we obtain
\[ \Theta^t_{\rho,u} = \epsilon(G,A_G)\theta(z) \sum_{\lambda \in [r^{-1}(u)]} \lambda(t) \sum_{Q}
   R(G^{u}_{\gamma_s},S_Q,1)(\gamma_u)[\phi_{Q_\lambda,Q}]_*\theta_\lambda(\gamma_s) \]
where the sum runs over the $G^{u}_{\gamma_s}(F)$-conjugacy classes inside the intersection of the $G^{u}(F)$-conjugacy class of $Q_\lambda$ with
$\tx{Lie}(G^{u}_{\gamma_s})(F)$. We obviously have
\[ [\phi_{Q_\lambda,Q}]_*\theta_\lambda = [\phi_{Q_\lambda,Q}]_*[\phi_{Q_0,Q_\lambda}]_*\theta_0 = [\phi_{Q_0,Q}]_*\theta_0 \]
and thus
\[ \Theta^t_{\rho,u} = \epsilon(G,A_G)\theta(z) \sum_{\lambda \in [r^{-1}(u)]} \lambda(t) \sum_{Q}
   R(G^{u}_{\gamma_s},S_Q,1)(\gamma_u)[\phi_{Q_0,Q}]_*\theta_0(\gamma_s) \]
Applying again Lemma \ref{lem:list} we obtain
\[ \Theta^t_{\rho,u} = \epsilon(G,A_G)\theta(z) \sum_{Q'} \langle \tx{inv}(Q_0,Q'),t_{q_0} \rangle^{-1} \sum_{Q}
   R(G^{u}_{\gamma_s},S_Q,1)(\gamma_u)[\phi_{Q_0,Q}]_*\theta_0(\gamma_s) \]
where now $Q'$ runs over a set of representatives for the $G^u(F)$-classes inside the $G^u$-stable class of $Q_0$, and $Q$ runs over a set of representatives
for the $G^u_{\gamma_s}(F)$-classes inside the intersection of the $G^u(F)$-class of $Q'$ with $\tx{Lie}(G^u_{\gamma_s})(F)$.

For any $Q'$ in the first summation set and $Q$ in the second, we have
\[ \tx{inv}(Q_0,Q') = \tx{inv}(Q_0,Q) \]
because $Q'$ and $Q$ are $G^u(F)$-conjugate. Thus we obtain
\[ \Theta^t_{\rho,u} = \epsilon(G,A_G)\theta(z) \sum_{Q} \langle \tx{inv}(Q_0,Q),t_{q_0} \rangle^{-1}
   R(G^{u}_{\gamma_s},S_Q,1)(\gamma_u)[\phi_{Q_0,Q}]_*\theta_0(\gamma_s) \]
where now $Q$ runs over a set of representatives for the $G^u_{\gamma_s}(F)$-classes inside the intersection of the stable class of $Q_0$ with
$\tx{Lie}(G^u_{\gamma_s})(F)$.

Now consider two elements $Q_1,Q_2$ in the above summation set, and assume that they are $G^u_{\gamma_s}$-conjugate. This means that
$\phi_{Q_1,Q_2}=\tx{Ad}(g)$ with $g \in G^u_{\gamma_s}$. Since $\gamma_s \in S_{Q_1}$ the expression $\phi_{Q_1,Q_2}(\gamma_s)$ is defined and we conclude that
it equals $\gamma_s$. Thus
\begin{eqnarray*}
[\phi_{Q_0,Q_1}]_*\theta_0(\gamma_s)&=&[\phi_{Q_2,Q_1}]_*[\phi_{Q_0,Q_2}]_*\theta_0(\gamma_s)\\
&=&[\phi_{Q_0,Q_2}]_*\theta_0( \phi_{Q_2,Q_1}^{-1}(\gamma_s) )\\
&=&[\phi_{Q_0,Q_2}]_*\theta_0(\gamma_s)
\end{eqnarray*}

Rearranging terms again we arrive at the desired formula for $\Theta^t_{\rho,u}$. \qed


\section{Character identities} \label{sec:chid}

In this section we assume all the notation established in the previous sections, in particular all parts of Section \ref{sec:main}. Our goal is to prove
Theorem \ref{thm:chid}.


\subsection{Beginning of the proof of \ref{thm:chid}}

\begin{lem} \label{lem:red0} Let $D$ be a diagonalizable group defined over $F$ and split over $F^u$. Then
\[ D(F) = D_s(F) \cdot \left[ D(F) \cap \bigcap_{\chi \in X^*(D)} [\tx{ker}(v\circ\chi)] \right] \]
where $D_s$ is the maximal split subtorus of $D$.
\end{lem}
\pf For any $x \in D(F^u)$ the map
\[ \lambda_x : X^*(D) \rightarrow \Z,\qquad \chi \mapsto v(\chi(x)) \]
is $\Z$-linear. Sending $x$ to $\lambda_x$ defines a homomorphism
\[ D(F^u) \rightarrow X_*(D) \]
which is $\Gamma$-equivariant because of the $\Gamma$-invariance of $v : [F^u]^\times \rightarrow \Z$. A right inverse of this homomorphism is given by
evaluation at $\pi$.

Now let $x \in D(F)$. Then $\lambda_x \in X_*(D)$ is $\Gamma$-fixed, and so its image $y=\pi^{\lambda_x} \in D(F^u)$ under the right inverse lies in $D_s(F)$.
Thus $x=y\cdot (xy^{-1})$ is the desired decomposition.\qed

\begin{lem} \label{lem:red1} Assume that $\gamma$ does not belong to $Z(F)G^u(F)_0$. Then both sides of \eqref{eq:chid} vanish. \end{lem}
\pf The left hand vanishes by \cite[Lemma 9.3.1]{DR09}. Turning to the right hand side, assume by way of contradiction that some $\gamma^H$ in the summation
set of \eqref{eq:chid} lies in $A_H(F)H(F)_0$, and write $\gamma^H = zx$. The admissible isomorphism $\phi_{\gamma^H,\gamma}$ sends $x$ into $G^u(F)_0$ and
because $H$ is elliptic it maps $A_H(F)$ to $A_G(F)$. Thus $\gamma=\phi_{\gamma^H,\gamma}(\gamma_H) \in A_G(F)G^u(F)_0$ contradicting the assumption of the
lemma. We conclude that all $\gamma^H$ occurring in the summation set of \eqref{eq:chid} lie outside of $A_H(F)H(F)_0$. But by the previous lemma,
$Z_H(F)H(F)_0=A_H(F)H(F)_0$,  because the set
\[ Z_H(F) \cap \bigcap_{\chi \in X^*(Z_H)} [\tx{ker}(v_F\circ\chi)] \]
lies in the maximal bounded subgroup of $T_0^H(F)$. By \cite[Lemma 9.3.1]{DR09} the left hand side of \eqref{eq:chid} also vanishes.\qed

\begin{lem} \label{lem:red2}
The isomorphism
\[ T_0^w \stackrel{\eta}{\longrightarrow} T_0^{H,w^H} \]
is defined over $F$. If $\gamma \in T_0^w(F)$ is topologically semi-simple and $z \in Z^\circ(F)$ then
\[ \theta(\gamma) = \theta^H(\eta(\gamma))\qquad \theta(z)=\theta^H(\eta(z))\lambda^G(z) \]
where $\lambda^G : Z^\circ(F) \rightarrow \C^\times$ is the character of \cite[Lemma 4.4.A]{LS87}.
\end{lem}
\pf Recall that $T_0^w$ is the torus whose complex dual is given by the complex torus $\hat T_0$ with $\Gamma$-action
\[ \sigma(t) = \tx{Ad}(\phi(\sigma))t \]
for all $\sigma \in W_F, t \in \hat T_0(\C)$ where the conjugation takes place in $^LG$. Analogously we have the torus $T_0^{H,w^H}$ whose complex dual is the
complex torus $\hat T_0^H$ with $\Gamma$-action
\[ \sigma(t) = \tx{Ad}(\phi^H(\sigma))t \]
for all $\sigma \in W_F, t \in \hat T_0^H(\C)$ where now the conjugation takes place in $^LH$. The statement that
\[ \eta : T_0^w \rightarrow T_0^{H,w^H} \]
is defined over $F$ is equivalent to the statement that the isomorphism of complex tori
\[ \hat\eta : \hat T_0^H \rightarrow \hat T_0 \]
is equivariant with respect to the above actions. But
\begin{eqnarray*}
\hat\eta(\tx{Ad}(\phi^H(\sigma))t)&=&{^L\eta}(\tx{Ad}(\phi^H(\sigma))t )\\
&=&\tx{Ad}({^L\eta}\phi^H(\sigma)){^L\eta}(t)\\
&=&\tx{Ad}(\phi(\sigma))\hat\eta(t)
\end{eqnarray*}

This proves the first assertion.

The restriction of $\theta$ to the maximal bounded subgroup of $T_0^\omega$, to which $\gamma$ belongs, is determined by the restriction of the Langlands
parameter $\phi_T$ of $\theta$ to inertia, which by construction is simply given by the restriction to inertia of $\phi={^L\eta}\circ\phi^H$. This restriction
is the cocycle
\begin{diagram}
I_F&\rTo^{\phi^H}&^LH&\rTo^{^L\eta}&^LG&\rTo&\hat G
\end{diagram}
which by construction lands in $\hat T_0$. Since $^L\eta$ is trivial on inertia, we see that this is the same as the cocycle
\begin{diagram}
I_F&\rTo^{\phi^H}&^LH&\rTo&\hat H&\rTo^{\hat\eta}&\hat G
\end{diagram}
which also lands in $\hat T_0$ and equals the restriction to inertia of $\hat\eta \circ \phi_{T^H}$. The latter is the cocycle determining the restriction of
$\theta^H \circ \eta$ to the maximal bounded subgroup of $T_0^\omega$. This proves the second assertion.

Let $T$ be any torus of $G$ coming from $H$. In \cite[\S3.5]{LS87} Langlands and Shelstad construct an element $a \in H^1(W_F,\hat T)$. The character
$\lambda^G(z)$ is then the restriction to $Z^\circ(F)$ of the character on $T(F)$ corresponding via the Langlands correspondence to $a$. The construction of
$a$ involves $\chi$-data, but one sees easily that its image under
\[ H^1(W_F,\hat T) \rightarrow H^1(W_F,\hat{Z^\circ}) \]
is independent of that choice and is in fact represented by the cocycle
\begin{diagram}
W_F&\rInto&^LH&\rTo^{^L\eta}&^LG&\rTo&^L[Z^\circ]&\rTo&\hat{Z^\circ}
\end{diagram}
By construction of the Langlands parameter $\phi_T$ of $\theta$, the restriction of $\theta$ to $Z^\circ(F)$ is given by the cocycle
\begin{diagram}
W_F&\rTo^{\phi^H}&^LH&\rTo^{^L\eta}&^LG&\rTo&^L[Z^\circ]&\rTo&\hat{Z^\circ}
\end{diagram}
while that of $\theta^H \circ \eta$ is given by the cocycle
\begin{diagram}
W_F&\rTo^{\phi^H}&^LH&\rTo&\hat H&\rTo^{\hat\eta}&\hat G&\rTo&\hat{Z^\circ}
\end{diagram}
It is clear that of these three cocycles, the second one equals the product of the first and the third, which implies the final statement of the lemma.\qed

\begin{cor} \label{cor:redcmp}
If equation \eqref{eq:chid} holds for all strongly regular semi-simple $\gamma \in G^u(F)_0$, then it holds for all strongly-regular semi-simple $\gamma \in
G^u(F)$.
\end{cor}
\pf By Lemma \ref{lem:red1} equation \eqref{eq:chid} holds trivially if $\gamma$ does not belong to $Z(F)G^u(F)_0$. By Lemma \ref{lem:red0} we have
$Z(F)G^u(F)_0=A_G(F)G^u(F)_0$, so it is enough to prove equation \eqref{eq:chid} for strongly regular semi-simple elements $\gamma=z\gamma'$ with $z \in
A_G(F)$ and $\gamma' \in G^u(F)_0$. By Proposition \ref{pro:schar} we know the behavior of the unstable character under central translations, namely
$\Theta^s_{\rho,u}(z\gamma) = \theta(z)\Theta^s_{\rho,u}(\gamma)$ and thus using our assumption we have
\begin{eqnarray*}
\Theta^s_{\rho,u}(z\gamma)&=&\theta(z)\sum_{\gamma^H \in H_\tx{sr}(F)/\tx{st}} \Delta_{\psi,u}(\gamma^H,\gamma)
    \frac{D(\gamma^H)^2}{D(\gamma)^2}\mc{S}\Theta_{\phi^H}(\gamma^H)\\
&=&\theta^H(\eta(z))\lambda^G(z)\sum_{\gamma^H \in H_\tx{sr}(F)/\tx{st}} \Delta_{\psi,u}(\gamma^H,\gamma)
    \frac{D(\gamma^H)^2}{D(\gamma)^2}\mc{S}\Theta_{\phi^H}(\gamma^H)
\end{eqnarray*}
where for the second equality we have invoked Lemma \ref{lem:red2}. Using \cite[Lemma 4.4.A]{LS87} and the obvious invariance of the terms $D(\gamma)$ and
$D^H(\gamma^H)$ under central translations this can be written as
\begin{eqnarray*}
&=&\sum_{\gamma^H \in H_\tx{sr}(F)/\tx{st}} \Delta_{\psi,u}(\eta(z)\gamma^H,z\gamma)
    \frac{D(\eta(z)\gamma^H)^2}{D(z\gamma)^2}\mc{S}\Theta_{\phi^H}(\eta(z)\gamma^H)\\
&=&\sum_{\gamma^H \in H_\tx{sr}(F)/\tx{st}} \Delta_{\psi,u}(\gamma^H,z\gamma)
    \frac{D(\gamma^H)^2}{D(z\gamma)^2}\mc{S}\Theta_{\phi^H}(\gamma^H)\\
\end{eqnarray*}
which was to be shown.\qed


\subsection{A reduction formula for the endoscopic lift of the stable character}

\begin{lem} \label{lem:d0}
Let $J$ be an unramified $F$-group and $y \in J(F)$ be a topologically semi-simple element belonging to a hyperspecial maximal compact subgroup.
Let $\gamma$ be an element of either $J(F)$ or $\tx{Lie}(J)(F)$ for which $\tx{Cent}(\gamma,J) \subset J_y$. Then the finite group $\pi_0(J^y(F))$ acts simply
transitively on the set of $J_y$-stable classes inside the $J^y$-stable class of $\gamma$.
\end{lem}
\pf Clearly $J^y(F)$ acts on the $J^y$-stable class of $\gamma$, and $J_y(F)$ preserves each $J_y$-stable class inside, so that we obtain an action of
$\pi_0(J^y)(F)$ on the set of $J_y$-stable classes inside the $J^y$-stable class of $\gamma$.

Consider the sequence
\[  1 \rightarrow J_y(F)\rightarrow J^y(F)\rightarrow\pi_0(J^y)(F)\rightarrow H^1(F,J_y)\rightarrow H^1(F,J^y) \]
By \cite[Prop 7.1]{Kot86} the last arrow has trivial kernel, which implies that the third arrow is surjective, so that we have a short exact sequence
\[ 1 \rightarrow J_y(F)\rightarrow J^y(F)\rightarrow\pi_0(J^y)(F) \rightarrow 1 \]

Let $\gamma'$ be $J^y$-stably conjugate to $\gamma$, and pick $g \in J^y(\ol{F})$ s.t. $\tx{Ad}(g)\gamma=\gamma'$ and $g^{-1}\sigma(g)\in J_\gamma \subset J_y$
for any $\sigma \in \Gamma$. Then the image $\bar g \in \pi_0(J^y)$ of $g$ belongs to $\pi_0(J^y)(F)$. Let $h \in J^y(F)$ be a preimage of $\bar g$. Then
$\tx{Ad}(h)\gamma$ and $\gamma'$ are stably conjugate by $gh^{-1} \in J_y(\ol{F})$. This proves transitivity.

To show simplicity, let $\gamma'$ by $J_y$-stably conjugate to $\gamma$ and pick $h \in J_y(\ol{F})$ s.t. $\tx{Ad}(h)\gamma=\gamma'$. If $g \in J^y(F)$ is also
an element s.t. $\gamma'=\tx{Ad}(g)\gamma$, then $gh^{-1} \in \tx{Cent}(\gamma,J) \subset J_y$ so $g \in J_y(\ol{F}) \cap J^y(F) = J_y(F)$. \qed

\rmk The same proof shows that under the same assumptions, $\pi_0(J^y(F))$ acts simply transitively on the set of $\tx{Ad}J_y(\ol{F})$-orbits in
$\tx{Ad}J^y(\ol{F})\gamma \cap J_y(F)$.

\begin{lem} \label{lem:d1}
Let $\gamma' \in G^u(F)$ be a strongly-regular semi-simple element. Assume that for some $\lambda \in r^{-1}(u)$ we have $\gamma'_s \in S_\lambda(F)$. Then
there exists a $\gamma \in G(F)$ stably conjugate to $\gamma'$ s.t. $\gamma_s \in S_0(F)$.
\end{lem}
\pf By construction we know that $\tx{Ad}(q_0q_\lambda^{-1}) : S_\lambda \rightarrow S_0$ is an isomorphism over $F$. Put $t =
\tx{Ad}(q_0q_\lambda^{-1})\gamma'_s$. Then $t$, being topologically semi-simple, belongs to the maximal bounded subgroup of $S_0(F)$ and thus lies in $G(O_F)$.
It follows form \cite[Prop. 7.1]{Kot86} that $G_t$ is quasi-split. The map $\tx{Ad}(q_0q_\lambda^{-1}) : G^u_{\gamma_s} \rightarrow G_t$ is a twist and so
there exists an $i \in G_t(\ol{F})$ s.t. if $T' = \tx{Cent}(\gamma',G_{\gamma'_s}^u)$ and $f=\tx{Ad}(iq_0q_\lambda^{-1})$ then the torus $T := f(T')$ and the
isomorphism $f : T' \rightarrow T$ are defined over $F$. Put $\gamma = f(\gamma')$. By construction $\gamma_s=t$ and $f$ is a $(\psi,u)$-equivalent twist, so
$\gamma$ is the element we want.\qed

\rmk The same proof can be applied to an element $\gamma^H \in H(F)$ and the trivial twist $(\tx{id},1) : H \rightarrow H$.

\begin{lem} \label{lem:dtf}
Let $\gamma \in G(F)_0$ and $\gamma^H \in H(F)_0$ be a pair of related strongly $G$-regular elements s.t. $\gamma_s \in S_0(F)$ and $\gamma^H_s \in S_0^H(F)$.
Then the admissible isomorphism $\phi_{\gamma^H,\gamma}$ makes $H_{\gamma^H_s}$ into an endoscopic group for $G_{\gamma_s}$. If $\Delta_0$ and $\Delta_0'$
denote the transfer factors for $(G,H)$ and $(G_{\gamma_s},H_{\gamma^H_s},\phi_{\gamma^H,\gamma})$ normalized with respect to admissible splittings (in the
sense of \cite[\S7]{Hal93}) then one has
\[ \Delta_0(\gamma^H,\gamma) = \Delta_0'(\gamma^H_u,\gamma_u) \]
\end{lem}
\pf By \cite[Prop. 7.1]{Kot86} both groups $H_{\gamma^H_s}$ and $G_{\gamma_s}$ are unramified, so they fall in the framework of \cite{Hal93} and one has the
normalization $\Delta_0'$ of the transfer factor with respect to an admissible splitting. We want to apply \cite[Thm. 10.18]{Hal93} to conclude
\[ \Delta_0(\gamma^H,\gamma) = \Delta_0'(\gamma^H,\gamma) \]
This theorem has two requirements. One is $p>e_G+1$, which is given in the statement of the theorem, and which we are assuming. The other one is $\gamma \in
G(O_F)$ and $\gamma^H \in H(O_F)$, which is a general requirement for the whole section 10 in loc. cit. However, tracing through the arguments of that section
one sees that up to the proof of Thm. 10.18, the only property of $\gamma$ and $\gamma^H$ which is used is that fact that they are compact and so have a
topological Jordan decomposition, together with the fact that homomorphisms preserve the topological Jordan decomposition and the knowledge of its explicit
form for elements of extensions of $F$. In the proof of Thm. 10.18 the elements $\gamma^H$ and $\gamma$ are replaced by high powers of themselves, let's call
them $\gamma'^H$ and $\gamma'$, which are very close (and can be made arbitrarily close) to $\gamma^H_s$ and $\gamma_s$. Then Lemma 8.1. of loc. cit. is
involved for the pair $(\gamma'^H,\gamma')$. For that Lemma it is essential that $\gamma'^H \in H(O_F)$ and $\gamma' \in G(O_F)$. But from our assumptions it
follows that $\gamma^H_s \in H(O_F)$ and $\gamma_s \in G(O_F)$, and since these groups are open, and the elements $\gamma'^H$ and $\gamma'$ can be made
arbitrarily close to $\gamma^H_s$ and $\gamma_s$, Lemma 8.1 can be invoked.

Thus we conclude that
\[ \Delta_0(\gamma^H,\gamma) = \Delta_0'(\gamma^H,\gamma) \]
By \cite[\S3.5]{LS90} there exists a character $\lambda : Z_{G_{\gamma_s}}(F) \rightarrow \C^\times$ s.t. for all strongly regular elements $z \in
H_{\gamma^H_s}(F)$ and $w \in G_{\gamma_s}(F)$ one has
\[ \Delta_0'(z\gamma^H_s,w\gamma_s) = \lambda(\gamma_s)\Delta_0'(z,w) \]
Thus
\[ \frac{\Delta_0'(\gamma^H,\gamma)}{\Delta_0'(\gamma^H_u,\gamma_u)} = \lambda(\gamma_s) = \frac{\Delta_0'(z\gamma^H_s,w\gamma_s)}{\Delta_0'(z,w)} \]

We choose $w$ to lie in an unramified torus $T \subset G_{\gamma_s}$. Then
\[ \frac{\Delta_0'(z\gamma^H_s,w\gamma_s)}{\Delta_0'(z,w)} = \langle a,\gamma_s \rangle \]
where $\langle a,\cdot \rangle$ is a character $T(F) \rightarrow \C^\times$ constructed in \cite[\S3.5]{LS87}. By \cite[Lemma 11.2]{Hal93} this character is
unramified, and thus takes trivial value at $\gamma_s$. It follows that
\[ \Delta_0'(\gamma^H,\gamma) = \Delta_0'(\gamma^H_u,\gamma_u) \]
and the proof is complete. \qed

\begin{lem} \label{lem:rhs}
For $\gamma \in G(F)_0$ the expression
\begin{equation} \label{ex:rhs}
\sum_{\gamma^H \in H_\tx{sr}(F)/\tx{st}} \Delta_0(\gamma^H,\gamma) \frac{D(\gamma^H)^2}{D(\gamma)^2} \mc{S}\Theta_{\phi^H}(\gamma^H)
\end{equation}
is equal to
\begin{equation} \label{ex:rhs1}
\sum_y \sum_\xi |\pi_0(H^y(F))|^{-1} \sum_{z \in H_y(F)_\tx{sr}/\tx{st}} \Delta_{0,y,\xi}(z,\gamma_u) \frac{D_{H_y}(z)^2}{D_{G_{\gamma_s}}(\gamma_u)^2}
    \mc{S}\Theta_{\phi^H}(yz) \end{equation}
where $y$ runs over a subset of $S_0^H(F)$ consisting of representatives for the stable classes of preimages of $\gamma_s$ which lie in $S_0^H(F)$, $\xi$ runs
over the $(H_y,G_{\gamma_s})$-equivalence classes of admissible embeddings mapping $y$ to $\gamma_s$, and $\Delta_{0,y,\xi}$ is the absolute transfer factor
for $(H_y,G_{\gamma_s},\xi)$ normalized with respect to an admissible splitting.
\end{lem}

\pf The sum of the first expression runs over the set of stable classes of strongly-regular semi-simple elements in $H(F)$. If $\gamma^H \in H(F)$ is
strongly-regular semi-simple, but $\gamma^H_s$ does not lie in a torus which is stably conjugate to $S_0^H$, then by Proposition \ref{pro:schar} we have
$\mc{S}\Theta_{\phi^H}(\gamma^H)=0$. Moreover if $\gamma^H_s$ is not a preimage of $\gamma_s$, then $\gamma^H$ is not a preimage of $\gamma$ and so
$\Delta_0(\gamma^H,\gamma)=0$. Thus if $\Gamma^H \subset H(F)$ is the set of strongly-regular semi-simple elements $\gamma^H$ for which $\gamma^H_s$ is a
preimage of $\gamma_s$ and lies in a torus stably conjugate to $S_0^H$, then we may restrict the summation in the first expression to $\Gamma^H/\tx{st}$. Let
$Y \subset S_0^H(F)$ be a set of representatives for the stable classes of those elements of $S_0^H(F)$ which occur as the topologically semi-simple parts of
elements in $\Gamma^H$. We claim that we have a surjective map
\[ p : \Gamma^H/\tx{st} \rightarrow Y \]
By Lemma \ref{lem:d1} and the remark thereafter every stable class $\mc{C} \subset \Gamma^H$ has a representative $\gamma^H$ s.t. $\gamma^H_s \in S_0^H(F)$. By
construction there exists $y \in Y$ stably conjugate to $\gamma^H_s$. By \cite[Prop. 7.1]{Kot86} there exists $h \in H(O_F)$ s.t. $\tx{Ad}(h)\gamma^H_s=y$. But
then $\tx{Ad}(h)\gamma^H \in \mc{C}$. We see that there are elements in  $\mc{C}$ whose topologically semi-simple parts lie in $Y$. If $\gamma^H,\gamma'^H \in
\mc{C}$ are two such elements, then their stable conjugacy implies the stable conjugacy of their topologically semi-simple parts, but by construction of $Y$
this means that their topologically semi-simple parts are actually equal. Thus we may define $p(\mc{C})$ by choosing any $\gamma^H \in \mc{C}$ with $\gamma^H_s
\in Y$ and sending it to $\gamma^H_s$.

Next we claim that for every  $y \in Y$ we have a surjective map
\[ [H_y(F)]_{(H,y)-\tx{sr},\tx{tu}}/\tx{st} \rightarrow p^{-1}(y),\qquad z \mapsto yz \]
whose fibers are torsors under $\pi_0(H^y(F))$. Here $[H_y(F)]_{(H,y)-\tx{sr},\tx{tu}}$ denotes the set of topologically unipotent elements $z \in H_y(F)$ for
which $yz$ is $H$-strongly regular, and $\tx{st}$ is stable conjugacy in $H_y$. It is immediate that this map is well-defined and surjective. We claim that
each fiber constitutes a single $H^y$-stable class. If $z,z'$ lie in the same fiber, then there exists $h \in H(\ol{F})$ s.t. $\tx{Ad}h(yz)=yz'$. But then
$\tx{Ad}h(y)=y$, so $h \in H^y(\ol{F})$, and $\tx{Ad}h(z)=z'$, which shows that $z,z'$ lie in the same $H^y$-stable class. Conversely if $z,z'$ lie in the same
$H^y$-stable class then they clearly lie in the same fiber. From Lemma \ref{lem:d0} it now follows that the fibers are torsors under $\pi_0(H^y(F))$.

We conclude that expression \eqref{ex:rhs} is equal to
\begin{equation} \label{ex:rhs1a}
\sum_{y\in Y} |\pi_0(H^y(F))|^{-1} \sum_{z \in [H_y(F)]_{(H,y)-\tx{sr},\tx{tu}}/\tx{st}} \Delta_0(yz,\gamma) \frac{D(yz)^2}{D(\gamma)^2}\mc{S}\Theta_{\phi^H}(yz)
\end{equation}

Consider $y,z$ contributing to the above expression. If $(yz,\gamma)$ is not a pair of $(G,H)$-related elements, then $\Delta_0(yz,\gamma)=0$. Now assume that
$(yz,\gamma)$ is a related pair. Then $(z,\gamma_u)$ is a pair of $(G_{\gamma_s},H_y,\phi_{yz,\gamma})$-related elements, and from Lemma \ref{lem:dtf} we know
that
\[ \Delta_0(yz,\gamma) = \Delta_{0,y,\phi_{yz,\gamma}}(z,\gamma_u) \]
Moreover, if $\xi$ is a $(G,H)$-admissible embedding carrying $y$ to $\gamma_s$ but not equivalent to $\phi_{yz,\gamma}$, the pair $(z,\gamma_u)$ is not
$(G_{\gamma_s},H_y,\xi)$-related, and thus $\Delta_{0,y,\xi}(z,\gamma_u)=0$. It follows that
\[ \Delta_0(yz,\gamma) = \sum_\xi\Delta_{0,y,\xi}(z,\gamma_u) \]
where $\xi$ runs over the set of $(G_{\gamma_s},H_y)$-equivalence classes of $(G,H)$-admissible embeddings carrying $y$ to $\gamma_s$. As in the proof of
\cite[Lem. 8.1]{Hal93} we have
\[ D(\gamma) = D_{G_{\gamma_s}}(\gamma_u)\qquad \tx{and} \qquad D(yz) = D_{H_y}(z) \]

Thus expression \eqref{ex:rhs1a} equals
\begin{equation*}
\sum_{y\in Y} \sum_\xi |\pi_0(H^y(F))|^{-1} \sum_{z \in [H_y(F)]_{(H,y)-\tx{sr},\tx{tu}}/\tx{st}} \Delta_{0,y,\xi}(z,\gamma_u) \frac{D_{H_y}(z)^2}{D_{G_{\gamma_s}}(\gamma_u)^2}\mc{S}\Theta_{\phi^H}(yz)
\end{equation*}

Finally, note that every $z \in H_y(F)_\tx{sr}$ which is a $(G_{\gamma_s},H_y,\xi)$-preimage of $\gamma_u$ automatically belongs to the set
$[H_y(F)]_{(H,y)-\tx{sr},\tx{tu}}$. Hence we may extend the summation over $z$ to all of $H_y(F)_\tx{sr}$. Also if $y \in S_0^H(F)$ is a preimage of $\gamma_s$
but does not belong to $Y$, then $H_y(F)$ does not contain a $(G_{\gamma_s},H_y,\xi)$-preimage of $\gamma_u$ for any $\xi$, and thus the terms
$\Delta_{0,y,\xi}(z,\gamma_u)$ vanish for all $\xi$ and $z$. Hence we may add to $Y$ representatives for the stable classes of such elements without changing
the value of the sum. But then the expression we obtain is \eqref{ex:rhs1}. \qed

\begin{cor} \label{cor:redcmp1}
If $\gamma \in G^u(F)_0$ is a strongly regular semi-simple element which does not have a stable conjugate $\gamma' \in G(F)_0$ with $\gamma'_s \in S_0(F)$,
then both sides of Equation \eqref{eq:chid} vanish.
\end{cor}
\pf Consider first the left hand side. In view of Proposition \ref{pro:schar}, it vanishes unless $\gamma_s$ lies in the centralizer of some $Q \in
\tx{Lie}(G^u)(F)$ stably conjugate to $Q_0$. But such a $Q$ is then rationally conjugate to $Q_\lambda$ for some $\lambda \in r^{-1}(u)$ and hence replacing
$\gamma$ by a rational conjugate we may assume $\gamma_s \in S_\lambda$. Thus, by Lemma \ref{lem:d1}, the non-vanishing of the left hand side of Equation
\eqref{eq:chid} implies the existence of $\gamma'$ as claimed.

We now turn to the right hand side. Let $\tilde\gamma \in G(F)_0$ be any stable conjugate of $\gamma$. By Lemma \ref{lem:rhs}, the right hand side of Equation
\eqref{eq:chid} vanishes at $\tilde\gamma$ unless there exists a triple $(y,\xi,z)$ s.t. $y \in S_0^H(F)$ is a preimage of $\tilde\gamma_s$, $\xi$ is a
$(G,H)$-admissible embedding s.t. $\xi(y)=\tilde\gamma_s$, and $z \in H_y(F)$ is a $(G_{\tilde\gamma_s},H_y,\xi)$-preimage of $\tilde\gamma_u$. By Lemma
\ref{lem:red2} the map
\begin{diagram}
S_0^H&\rTo^{\tx{Ad}(q_0^H)^{-1}}&T_0^{H,w^H}&\rTo^{\eta^{-1}}&T_0^w&\rTo^{\tx{Ad}(q_0)}&S_0
\end{diagram}
is an admissible isomorphism defined over $F$. Let $y'$ be the image of $y$ under this isomorphism. Then $G_{y'}$ is quasi-split by \cite[Prop. 7.1]{Kot86} and
so there exists $z' \in G_{y'}(F)$ which is an image of $z$. But then $\gamma'=y'z'$ is a stable conjugate of $\tilde\gamma$, hence of $\gamma$. Thus the
non-vanishing of the right hand side of Equation \eqref{eq:chid} at $\tilde\gamma$ implies the existence of $\gamma'$ as claimed. But since $\gamma$ and
$\tilde\gamma$ are stably conjugate, the non-vanishing of said expression at $\gamma$ is equivalent to its non-vanishing at $\tilde\gamma$, since the value at
$\tilde\gamma$ differs from the value at $\gamma$ by a non-zero multiplicative factor.\qed


\subsection{Lemmas about transfer factors}

In this section $G'$ is an unramified $F$-group and $(H',s,{^L\eta})$ is an unramified extended endoscopic triple for $G'$. Let $(T'_0,B'_0)$ be a Borel pair
of $G'$ over $F$. We choose hyperspecial points in the buildings of $G'$ and $H'$, s.t. the one for $G'$ lies in the apartment of $T'_0$. We also choose an
admissible splitting $(T'_0,B'_0,\{X'_\alpha\})$ for $G'$ in the sense of \cite[\S7]{Hal93}. Then we have the transfer factors normalized with respect to that
splitting both on the group level (\cite[\S3.7]{LS87}), as well as on the Lie algebra level (\cite{Kot99}). We will call both these transfer factors
$\Delta_0$, as there will be no possibility of confusion between the two.

\begin{lem} \label{lem:top1}
For any semi-simple strongly regular topologically unipotent $\gamma^H \in H'(F)$ and $\gamma \in G'(F)$, we have
\[ \Delta_0(\gamma^H,\gamma) \frac{D(\gamma^H)^2}{D(\gamma)^2} = \Delta_0(\log(\gamma^H),\log(\gamma))\frac{D(\log(\gamma^H))}{D(\log(\gamma))} \]
\end{lem}
\pf We choose a positive integer $m$ with the property that the sequences
\[ \gamma_k = \gamma^{p^{km}}\qquad \gamma^H_k = [\gamma^H]^{p^{km}} \]
converge to $1$ (cf. \cite[\S7]{DR09}), and put
\[X^H=\log(\gamma^H), X=\log(\gamma),X^H_k = \log(\gamma^H_k), X_k=\log(\gamma_k)\]

As argued in \cite[\S10]{Hal93} we have
\[\Delta_0(\gamma^H_{2k},\gamma_{2k}) \frac{D(\gamma^H_{2k})^2}{D(\gamma_{2k})^2} = |p^{km}|^{-N}\Delta_0(\gamma^H,\gamma) \]
where $N$ is the number of roots in $G'$ outside $H'$ and $|\ |$ is the unique absolute value on $\ol{F}$ extending that on $F$. By the same arguments one also
has
\[ \Delta_0(X_{2k}^H,X_{2k})\frac{D(X_{2k}^H)}{D(X_{2k})} = |p^{km}|^{-N}\Delta_0(X^H,X)\frac{D(X^H)}{D(X)} \]
Thus it will be enough to show the equality claimed in the lemma with $\gamma^H,\gamma$ replaced by $\gamma^H_{2k}, \gamma_{2k}$ for some $k$ which we may
freely choose.

As argued in \cite[\S2.3]{Wal97}, there exists a positive integer $K$ s.t. for all $k>K$
\[ \Delta_0(\gamma^H_{2k},\gamma_{2k})\frac{D(\gamma^H_{2k})}{D(\gamma_{2k})} = \Delta_0(X^H_{2k},X_{2k}) \]
We now claim that, after potentially increasing $K$, we have
\[ D(\gamma^H_{2k}) = D(X^H_{2k}) \qquad D(\gamma_{2k}) = D(X_{2k}) \]
For this it is enough to show that if $T \subset G'$ is a maximal torus with Lie algebra $\mf{t} \subset \mf{g}'$ and $Y \in \mf{t}(F)$ is small enough then
for all roots $\alpha \in R(T,G')$ we have
\[ |\alpha(\exp(Y))-1| = |d\alpha(Y)| \]
Let $E/F$ be the extension splitting $T$, and let $v_E$ be the unique valuation $E$ extending that on $F$ (here we deviate from our usual notation). Then
\[ |\alpha(\exp(Y))-1| = |\exp(d\alpha(Y))-1| = |d\alpha(Y)+\sum_{k>1} \frac{d\alpha(Y)^k}{k!}| \]
Putting $u=d\alpha(Y)$, we have by a computation similar to the proof of \cite[B.1.1]{DR09}
\[ v_E(\frac{u^k}{k!}) = kv_E(u)-eA(k) > kv_E(u)-(k-1) \]
where as in loc. cit. $A(k) = \sum_{i>0} \lfloor\frac{k}{p^i}\rfloor$ and $e$ is the ramification degree of $F/\Q_p$. Thus if $v_E(u) \geq 1$ then for all
$k>1$
\[ v_E(\frac{u^k}{k!}) > v_E(u) \]
from which follows
\[ |u+\sum_{k>1} \frac{u^k}{k!}| = |u| \]
This finishes the proof of the claim about $D$ and the lemma follows.\qed

\begin{lem} \label{lem:top2}
Let $S \subset G'$ be a torus (as usual defined over $F$) which is defined and split over $O_{F^u}$. Let $Q \in \tx{Lie}(S)(O_{F^u}) \cap \mf{g}'(F)$ be
semi-simple regular, and $Q^H$ be any preimage of $Q$ in $\mf{h}'(F)$. Then
\[ \Delta_0(Q^H,Q) = 1 \]
\end{lem}
\pf For $\alpha \in R(S,G')$ let $a_\alpha = d\alpha(Q)$. As Kottwitz observes in \cite{Kot99}, this defines a-data for $R(S,G')$ and with respect to that
a-data, $\Delta_{II}(Q^H,Q)=1$. To show that $\Delta_I(Q^H,Q)=1$ we adapt the argument of \cite[Lem. 7.2]{Hal93}. Since $S$ is defined and split over
$O_{F^u}$, all roots of $S$ are defined over $O_{F^u}$ and we have $a_\alpha \in O_{F^u}$. Then as in loc. cit. we see that the cocycle $m(\sigma_{S})$
constructed in \cite[\S2.3]{LS87} takes values in $G'(O_{F^u})$. Since the torus $T'_0$ is also defined over $O_{F^u}$, there exists $g \in G'(O_{F^u})$ s.t.
$S=\tx{Ad}(g)T'_0$. Thus the cocycle $\tx{Ad}(g)^{-1}m(\sigma_{S})$ of $\Gamma$ in $S(\ol{F})$ takes values in $S(O_{F^u})$ and is thus cohomologically
trivial. But $\Delta_I(Q^H,Q)$ is the value of a character on $H^1(\Gamma,S)$ at that cocycle.\qed

\begin{lem} \label{lem:top3}
Let $\gamma^H \in H'(F)$ and $\gamma \in G'(F)$ be semi-simple, strongly regular, and topologically unipotent. Then $\gamma$ is an image of
$\gamma^H$ if and only if $\log(\gamma)$ is an image of $\log(\gamma^H)$.
\end{lem}
\pf We define $\gamma^H_k,\gamma_k,X^H,X,X^H_k,X_k$ as in the proof of Lemma \ref{lem:top1}. It is clear that $\gamma$ is an image of $\gamma^H$ if and only if
$\gamma_k$ is an image of $\gamma^H_k$ for some (then any) $k$. The same holds for the $X$'s. This reduces the proof to the case where the elements are near
the identity, in which case it is clear.\qed


\subsection{Completion of the proof of theorem \ref{thm:chid}}

By Corollaries \ref{cor:redcmp} and \ref{cor:redcmp1} it is enough to prove Equation \eqref{eq:chid} for all strongly regular semi-simple elements $\gamma \in
G^u(F)_0$ which have a stable conjugate $\gamma' \in G(F)_0$ s.t. $\gamma'_s \in S_0(F)$. Fix such a pair $\gamma,\gamma'$ and consider the value at $\gamma$
of the right hand side of Equation \eqref{eq:chid}:

\begin{equation}
\sum_{\gamma^H \in H_\tx{sr}(F)/\tx{st}} \Delta_{\psi,u}(\gamma^H,\gamma) \frac{D(\gamma^H)^2}{D(\gamma)^2}\mc{S}\Theta_{\phi^H}(\gamma^H) \label{ex:rhs10}
\end{equation}

By construction of $\Delta_{\psi,u}$ we have
\[ \Delta_{\psi,u}(\gamma^H,\gamma) = \epsilon_L(V,\psi)\Delta_0(\gamma^H,\gamma') \langle \tx{inv}(\gamma',\gamma), \hat\phi_{\gamma',\gamma^H}(s)
\rangle^{-1} \]
where $\Delta_0$ is the absolute transfer factor for $(G,H)$ normalized with respect to our chosen splitting. By Lemma \ref{lem:red2} the map
\begin{diagram}
S_0^H&\rTo^{\tx{Ad}(q_0^H)^{-1}}&T_0^{H,w^H}&\rTo^{\eta^{-1}}&T_0^w&\rTo^{\tx{Ad}(q_0)}&S_0
\end{diagram}
is an admissible isomorphism defined over $F$. We fix $Q_0 \in \tx{Lie}(S_0)(F)$ satisfying the requirements of the element $X_S$ in \cite[Lemma 12.4.3]{DR09},
and let $Q_0^H$ be the preimage of $Q_0$ under this embedding. Then $Q_0^H$ also satisfies the same requirements.

We now apply Lemma \ref{lem:rhs} and Proposition \ref{pro:schar} to conclude that \eqref{ex:rhs10} equals

\begin{eqnarray} \label{ex:rhs11}
&&\epsilon_L(V,\psi)\epsilon(H,A_H)\sum_y \sum_\xi |\pi_0(H^y(F))|^{-1} \sum_{P^H} [\phi_{Q^H_0,P^H}]_*\theta_0^H(y) \nonumber\\
&&\langle \tx{inv}(\gamma',\gamma), \hat\phi_{\gamma',\gamma^H}(s) \rangle^{-1} \sum_{z \in H_y(F)_\tx{sr}/\tx{st}} \Delta_{0,y,\xi}(z,\gamma'_u) \frac{D_{H_y}(z)^2}{D_{G_{\gamma'_s}}(\gamma'_u)^2} \nonumber\\
&&\sum_{Q^H} R(H_{y},S_{Q^H},1)(z)
\end{eqnarray}
Let us recall the summation sets. $y$ runs over a set $Y \subset S_0^H(F)$ representing the stable classes of preimages of $\gamma'_s$ which intersect
$S_0^H(F)$, $\xi$ runs over the $(G_{\gamma'_s},H_y$)-equivalence classes of $(G,H)$-admissible embeddings which map $y$ to $\gamma'_s$, $P^H$ runs over a set
of representatives for the $H_y$-stable classes of elements of $\tx{Lie}(H_y)(F)$ which are $H$-stably conjugate to $Q_0^H$, $z$ runs over the stable classes
of strongly regular elements in $H_y(F)$, and $Q^H$ runs over a set of representatives for the $H_y(F)$-classes inside the $H_y$-stable class of $P^H$.

Consider a triple $(y,\xi,P^H)$. Since $G_{\gamma'_s}$ is quasi-split, there exists an $(G_{\gamma'_s},H_y,\xi)$-image $P' \in \tx{Lie}(G_{\gamma'_s})(F)$ of
$P^H$, unique up to stable conjugacy. We claim that the map
\[ p : (y,\xi,P^H) \mapsto P' \]
is a surjection from the set of triples $(y,\xi,P^H)$ occurring in \eqref{ex:rhs11} to the set of $G_{\gamma'_s}$-stable classes of elements of
$\tx{Lie}(G_{\gamma'_s})(F)$ stably conjugate to $Q_0$, and moreover that the fiber of this surjection through $(y,\xi,P^H)$ is a torsor under $\pi_0(H^y(F))$
for the action of this group by conjugation on all factors of the triple (the first factor is of course fixed by this action).

To see surjectivity, choose $P'$ in the target of $p$. Let $\tilde y=\phi_{P',Q_0^H}(\gamma'_s)$. There exists a $y \in Y$ stably conjugate to $\tilde y$. By
\cite[Prop. 7.1]{Kot86} there exists $h \in H(O_F)$ s.t. $\tx{Ad}(h)\tilde y=y$. Put $P^H = \tx{Ad}(h)Q_0^H$. Then $(y,\phi_{P^H,P'},P^H)$ is a preimage of
$P'$ under $p$.

Now let $(y,\xi,P^H)$ be an element in the source of $p$ and let $P'$ be its image. We claim that the map
\[ \tilde p: \tilde P^H \mapsto (y,\phi_{\tilde P^H,P'},\tilde P^H) \]
is an $\pi_0(H^y(F))$-equivariant bijection from the set of $H_y$-stable classes inside the $H^y$-stable class of $P^H$ to the fiber of $p$ through
$(y,\xi,P^H)$. Once this has been shown, the claim about the fibers of $p$ will follow from Lemma \ref{lem:d0}.

Indeed, let $\tilde P^H$ be $H^y$-stably conjugate to $P^H$. Then $\phi_{\tilde P^H,P^H}(y)=y$ and moreover since $P'$ is a $(G_{\gamma'_s},H_y,\xi)$-image of
$P^H$ we have $\phi_{P^H,P'}(y)=\gamma'_s$. This implies $\phi_{\tilde P^H,P'}(y) = \gamma'_s$ and we see that $(y,\phi_{\tilde P^H,P'},\tilde P^H)$ belongs to
the target of the proposed map $\tilde p$. If $\tilde P^H$ is replaced by an $H_y$-stable conjugate, then $\phi_{\tilde P^H,P'}$ remains within its equivalence
class. We see that $\tilde p$ is a well-defined and $\pi_0(H^y(F))$-equivariant map as claimed. It is clearly injective. To show surjectivity, let $(\tilde
y,\tilde \xi,\tilde P^H) \in p^{-1}(P')$. By definition of the map $p$, we must have that $\tilde\xi$ and $\phi_{\tilde P^H,P'}$ are $(G_{\gamma'_s},H_{\tilde
y})$-equivalent and $\tilde y=\phi_{P',\tilde P^H}(\gamma'_s)$ and so we only have to show that $\tilde P^H$ and $P^H$ are $H^y$-stably conjugate. We have
$\phi_{P^H,P'}(y)=\gamma'_s = \phi_{\tilde P^H,P'}(\tilde y)$. But recall that $P^H$ and $\tilde P^H$ are $H$-stably conjugate. Thus $\phi_{P^H,\tilde P^H}$ is
defined and since $\phi_{\tilde P^H,P'}=\phi_{P^H,P'}\circ\phi_{\tilde P^H,P^H}$ we have $\phi_{P^H,\tilde P^H}(y)=\tilde y$. But $Y$ contains only one element
per stable class, which forces $y=\tilde y$, and so $\phi_{P^H,\tilde P^H}(y)=y$, i.e. $P^H$ and $\tilde P^H$ are $H^y$-stably conjugate. This conclude the
proof of the claim about the map $p$.

Consider a triple $(y,\xi,P^H)$ contributing to \eqref{ex:rhs11} and let $P'$ be its image under $p$. We focus on the part of \eqref{ex:rhs11} given by
\begin{eqnarray} \label{ex:rhs11p1}
&&\langle \tx{inv}(\gamma',\gamma), \hat\phi_{\gamma',\gamma^H}(s) \rangle^{-1} \sum_{z \in H_y(F)_\tx{sr}/\tx{st}} \Delta_{0,y,\xi}(z,\gamma'_u) \frac{D_{H_y}(z)^2}{D_{G_{\gamma'_s}}(\gamma'_u)^2} \nonumber\\
&&\sum_{Q^H} R(H_{y},S_{Q^H},1)(z)
\end{eqnarray}
The map $\phi_{\gamma',\gamma}$ defines an inner twist $G_{\gamma'_s} \rightarrow G^u_{\gamma_s}$ and maps $\gamma'_u$ to $\gamma_u$. From this it follows that
$D_{G_{\gamma'_s}(\gamma'_u)} = D_{G^u_{\gamma_s}(\gamma_u)}$, and $\tx{inv}(\gamma',\gamma)=\tx{inv}(\gamma'_u,\gamma_u)=\tx{inv}(X',X)$, where
$X'=\log(\gamma'_u)$, $X=\log(\gamma_u)$. All $z$ which are preimages of $\gamma'_u$ are topologically unipotent, so we may restrict the sum over $z$ to the
topologically unipotent elements. Put $Z=\log(z)$. We will use Lemma \ref{lem:top1} with $G'=G_{\gamma'_s}$ and $H'=H_y$. By \cite[Prop. 7.1]{Kot86} these
groups are unramified and come with fixed hyperspecial maximal compact subgroups. We see that \eqref{ex:rhs11p1} equals
\begin{eqnarray} \label{ex:rhs11p2}
&&\langle \tx{inv}(X',X), \hat\phi_{\gamma',\gamma^H}(s) \rangle^{-1} \sum_{Z \in \mf{h}_y(F)_\tx{sr}/\tx{st}} \Delta_{0,y,\xi}(Z,X') \frac{D_{\mf{h}_y}(Z)}{D_{\mf{g}_{\gamma'_s}}(X')} \nonumber \\
&&\sum_{Q^H} R(H_{y},S_{Q^H},1)(z)
\end{eqnarray}
The function
\[ \Delta_{0,y,\phi_{\gamma',\gamma}\circ\xi}(Z,X) := \Delta_{0,y,\xi}(Z,X') \langle \tx{inv}(X',X), \hat\phi_{\gamma',\gamma^H}(s) \rangle^{-1} \]
is a transfer factor for $(\mf{g}^u_{\gamma_s},\mf{h}_y,\phi_{\gamma',\gamma}\circ\xi)$. Applying \cite[Lem. 12.4.3]{DR09} we conclude that \eqref{ex:rhs11p2}
equals
\begin{eqnarray} \label{ex:rhs11p3}
\sum_{Z \in \mf{h}_y(F)_\tx{sr}/\tx{st}} \Delta_{0,y,\phi_{\gamma',\gamma}\circ\xi}(Z,X) \frac{D_{\mf{h}_y}(Z)}{D_{\mf{g}_{\gamma_s}}(X)} \sum_{Q^H} \epsilon(H_y,A_{H_y})\hat\mu^{H_y}_{Q^H}(Z)
\end{eqnarray}
Here $\hat\mu_{Q^H}^{H_y}$ is the Fourier transform (with respect to the transfer to $\mf{h}_y$ of the bilinear form $B$ and the character $\psi$) of the
orbital integral at $Q^H$ on $\mf{h}_y(F)$.

We will now apply \cite[Conj. 1.2]{Wal97}, which is now a theorem due to the work of \cite{Wal97}, \cite{Wal06}, \cite{HCL07} and \cite{Ngo08}. According to
it, \eqref{ex:rhs11p3} equals
\begin{eqnarray} \label{ex:rhs11p4}
\gamma_\psi(\mf{g}^u_{\gamma_s})\gamma_\psi(\mf{h}_y)^{-1}\epsilon(H_y,A_{H_y})\sum_Q \Delta_{0,y,\phi_{\gamma',\gamma}\circ\xi}(P^H,Q)\hat\mu^{G_{\gamma_s}}_Q(X)
\end{eqnarray}
where $Q$ runs over a set of representatives for the conjugacy classes of regular semi-simple elements in $\mf{g}_{\gamma_s}(F)$.

For a moment we consider the signs in \eqref{ex:rhs11p4}. The group $H_y$ contains $S_0^H$, which is an elliptic maximal torus of $H$. Thus the inclusion $Z_H
\rightarrow Z_{H_y}$ restricts to an isomorphism $A_H \rightarrow A_{H_y}$. The group $G_{\gamma'_s}$ contains $S_0$, which is an elliptic maximal torus of
$G$, and again we get an isomorphism $A_G \rightarrow A_{G_{\gamma'_s}}$. The group $G^u_{\gamma_s}$ is an inner twist of $G_{\gamma'_s}$ and so we have an
isomorphism $A_{G_{\gamma'_s}} \rightarrow A_{G^u_{\gamma_s}}$. Finally since $H$ is elliptic for $G$, the natural inclusion $Z_G \rightarrow Z_H$ restricts to
an isomorphism $A_G \rightarrow A_H$. All in all this gives an isomorphism $A_{H_y} \rightarrow A_{G^u_{\gamma_s}}$. Using this and the transitivity of the
sign $\epsilon(\cdot,\cdot)$ we conclude
\[\epsilon(H_y,A_{H_y}) = \epsilon(H_y,G_{\gamma'_s})\epsilon(G_{\gamma'_s},G^u_{\gamma_s})\epsilon(G^u_{\gamma_s},A_{G^u_{\gamma_s}}) \]
From \cite[\S12.3]{DR09} we know
\[\epsilon(G_{\gamma'_s},G^u_{\gamma_s}) = \gamma_\psi(\mf{g}_{\gamma'_s})\gamma_\psi(\mf{g}^u_{\gamma_s})^{-1} \]
while from Proposition \ref{pro:signs} we know
\[\epsilon(H_y,G_{\gamma'_s}) = \gamma_\psi(\mf{h}_y)\gamma_\psi(\mf{g}_{\gamma'_s})^{-1} \]
It follows that \eqref{ex:rhs11p4} equals
\begin{eqnarray} \label{ex:rhs11p5}
\sum_Q \Delta_{0,y,\phi_{\gamma',\gamma}\circ\xi}(P^H,Q) \epsilon(G^u_{\gamma_s},A_{G^u_{\gamma_s}})\hat\mu^{G_{\gamma_s}}_Q(X)
\end{eqnarray}
where $Q$ runs over the same set as in \eqref{ex:rhs11p4}.

Now there is a natural injection from the set of $G^u_{\gamma_s}$-stable classes of regular semi-simple elements in $\mf{g}^u_{\gamma_s}(F)$ stably conjugate
to $Q_0$ to the set of $G_{\gamma'_s}$-stable classes of regular semi-simple elements in $\mf{g}_{\gamma'_s}(F)$ stably conjugate to $Q_0$. If $P'$ is not in
the image of that injection, then \eqref{ex:rhs11p5} is zero. Otherwise let $P \in\mf{g}_{\gamma_s}(F)$ be an element whose class maps to that of $P'$. Then
\eqref{ex:rhs11p5} equals
\begin{eqnarray} \label{ex:rhs11p6}
\sum_{Q} \Delta_{0,y,\xi}(P^H,Q_0)\langle \tx{inv}(Q_0,Q),s_{q_0} \rangle^{-1} \epsilon(G^u_{\gamma_s},A_{G^u_{\gamma_s}})\hat\mu^{G_{\gamma_s}}_Q(X)
\end{eqnarray}
where $Q$ runs over the set of $G^u_{\gamma_s}(F)$-conjugacy classes inside the $G^u_{\gamma_s}$-stable class of $P$.

The torus $S_0 \subset G_{\gamma'_s}$ and the element $Q_0$ satisfy the requirements of Lemma \ref{lem:top2}. Moreover the element $Q$ satisfies the
requirements of \cite[Lem 12.4.3]{DR09} on the element $X_S$. Thus \eqref{ex:rhs11p6} equals
\begin{eqnarray} \label{ex:rhs11p7}
\sum_{Q} \langle \tx{inv}(Q_0,Q),s_{q_0} \rangle^{-1} R(G^u_{\gamma_s},S_Q,1)(\gamma_u)
\end{eqnarray}
where $Q$ runs over the same set as in \eqref{ex:rhs11p6}.

To recapitulate, for a triple $(y,\xi,P^H)$ contributing to $\eqref{ex:rhs11}$ there are two possibilities. Either $P'=p(y,\xi,P^H)$ does not lie in the image
of the natural injection from the regular semi-simple stable classes in $\mf{g}^u_{\gamma_s}(F)$ to those in $\mf{g}_{\gamma'_s}(F)$ given by the inner twist
$\phi_{\gamma',\gamma} : G_{\gamma'_s} \rightarrow G^u_{\gamma_s}$, in which case the summand corresponding to that triple is zero. Or it does lie in that
image, and if $P$ is an element of the stable class in $\mf{g}^u_{\gamma_s}(F)$ which maps to that of $P'$, then the summand of \eqref{ex:rhs11} corresponding
to $(y,\xi,P^H)$ equals \eqref{ex:rhs11p7}.

After restricting the sums in \eqref{ex:rhs11} to the subset of triples $(y,\xi,P^H)$ whose image under $p$ lies in the image of the natural injection of
stable classes provided by $\phi_{\gamma',\gamma}$, we obtain a map
\[ (y,\xi,P^H) \mapsto P \]
which is a surjection on the set of $G^u_{\gamma_s}$-stable classes of elements of $\mf{g}^u_{\gamma_s}(F)$ which are stably conjugate to $Q_0$, and the fiber
of that surjection through $(y,\xi,P^H)$ is a torsor under $\pi_0(H^y(F))$. This of course follows from the corresponding property of the map $p$.

Before we apply this to the expression \eqref{ex:rhs11}, we need to note that if $(y,\xi,P^H)$ maps to $P$, then since $\phi_{P,P^H}(\gamma_s)=y$ we have
\begin{eqnarray*}
[\phi_{Q_0^H,P^H}]_*\theta_0^H(y)&=&[\phi_{P^H,P}]_*[\phi_{Q_0^H,P^H}]_*\theta_0^H(\gamma_s)\\
&=&[\phi_{Q_0,P}]_*[\phi_{Q_0^H,Q_0}]_*\theta_0^H(\gamma_s) \\
&=&[\phi_{Q_0,P}]_*\theta_0(\gamma_s)
\end{eqnarray*}
where the last equality follows from Lemma \ref{lem:red2}.

With this in mind, we see that \eqref{ex:rhs11} equals
\begin{eqnarray} \label{ex:rhs12}
\epsilon_L(V,\psi)\epsilon(H,A_H)&\cdot&\sum_P [\phi_{Q_0,P}]_*\theta_0(\gamma_s) \nonumber \\
&&\sum_{Q} \langle \tx{inv}(Q_0,Q),s_{q_0} \rangle^{-1} R(G^u_{\gamma_s},S_Q,1)(\gamma_u)
\end{eqnarray}
where $P$ runs over a set of representatives for the $G^u_{\gamma_s}$-stable classes of elements in $\mf{g}^u_{\gamma_s}(F)$ which are $G^u$-stably conjugate
to $Q_0$, and $Q$ runs over a set of representatives for the $G^u_{\gamma_s}(F)$-conjugacy classes inside the $G^u_{\gamma_s}$-stable class of $P$.

Again using the transitivity of $\epsilon(\cdot,\cdot)$ and the isomorphism $A_G\cong A_H$ we can write
\[ \epsilon(H,A_H)=\epsilon(H,G)\epsilon(G,A_G) \]
and thus using Proposition \ref{pro:signs} we see that \eqref{ex:rhs12} equals
\begin{eqnarray} \label{ex:rhs13}
\epsilon(G,A_G)\sum_P [\phi_{Q_0,P}]_*\theta_0(\gamma_s) \sum_{Q} \langle \tx{inv}(Q_0,Q),s_{q_0} \rangle^{-1} R(G^u_{\gamma_s},S_Q,1)(\gamma_u)
\end{eqnarray}
with both sums as in \eqref{ex:rhs12}. By Proposition \ref{pro:schar} this is the left hand side of Equation \eqref{eq:chid}. This completes the proof of
Theorem \ref{thm:chid}. \qed


\newpage

\begin{thebibliography}{9999999}
\setlength{\parindent}{0pt} \setlength{\parskip}{0pt} \footnotesize

\bibitem[Deb06]{Deb06}  S. DeBacker, Parameterizing conjugacy classes of maximal unramified tori via Bruhat-Tits theory.  Michigan Math. J. 54 (2006), no. 1,
    157--178
\bibitem[DR09]{DR09}    S. DeBacker, M. Reeder, Depth-zero supercuspidal $L$-packets and their stability.  Ann. of Math. (2)  169  (2009),  no. 3, 795--901
\bibitem[Hal93]{Hal93}  T. C. Hales, A simple definition of transfer factors for unramified groups, Contemporary Math., 145, (1993) 109--134
\bibitem[HCL07]{HCL07}  T. C. Hales, R. Cluckers, F. Loeser, Transfer Principle for the Fundamental Lemma, preprint, arXiv:0712.0708
\bibitem[Kot83]{Kot83}  R. E. Kottwitz, Sign changes in harmonic analysis on reductive groups.  Trans. Amer. Math. Soc.  278  (1983), no. 1, 289--297
\bibitem[Kot86]{Kot86}  R. E. Kottwitz, Stable trace formula: Elliptic singular terms, Math. Ann. 275 (1986), 365--399
\bibitem[Kot99]{Kot99}  R. E. Kottwitz, Transfer factors for Lie algebras.  Represent. Theory  3  (1999), 127--138
\bibitem[KS99]{KS99}    R. E. Kottwitz, D. Shelstad, Foundations of twisted endoscopy.  Astérisque  No. 255  (1999)
\bibitem[KV1]{KV1}      D. Kazhdan, Y. Varshavsky, Endoscopic decomposition of certain depth zero representations.  Studies in Lie theory,  223--301,
    Progr. Math., 243, Birkhäuser Boston, Boston, MA, 2006
\bibitem[KV2]{KV2}      D. Kazhdan, Y. Varshavsky, On endoscopic transfer of Deligne-Lusztig functions, preprint, arXiv:0902.3426
\bibitem[LS87]{LS87}    R. P. Langlands, D. Shelstad, On the definition of transfer factors, Math. Ann., vol. 278 (1987), 219--271
\bibitem[LS90]{LS90}    R. P. Langlands, D. Shelstad, Descent for transfer factors, in the Grothendieck Festschrift, Vol. II, Birkh\"auser (1990) pp. 485--563.
\bibitem[Ngo08]{Ngo08}  B. C. Ngo, Le lemme fondamental pour les algebres de Lie, preprint, arXiv:0801.0446v3
\bibitem[Ser79]{Ser79}  J. P. Serre, Local fields, Springer-Verlag, 1979.
\bibitem[Tat77]{Tat77}  J. Tate, Number theoretic background, Automorphic Forms, Representations and $L$-Functions, part 2, Proc. Sympos. Pure Math., XXXIII,
    Amer. Math. Soc., Providence, R.I., 1977, pp. 3--26.
\bibitem[Wal95]{Wal95}  J.-L. Waldspurger, Une formule des traces locale pour les algèbres de Lie $p$-adiques, J. Reine Angew. Math.  465  (1995), 41--99.
\bibitem[Wal97]{Wal97}  J.-L. Waldspurger, Le lemme fondamental implique le transfert, Compositio Math. 105 (1997), 153--236.
\bibitem[Wal06]{Wal06}  J.-L. Waldspurger, Endoscopie et changement de caractéristique, J. Inst. Math. Jussieu  5  (2006),  no. 3, 423--525.
\end{thebibliography}
\end{document}